\documentclass[letterpaper,10pt]{article}

\usepackage[T1]{fontenc}
\usepackage[utf8]{inputenc}

\usepackage[symbol]{footmisc}

\setfnsymbol{wiley}

\usepackage{mathrsfs}

\usepackage{framed}

\usepackage[titletoc]{appendix}

\usepackage{xspace}

\usepackage{etex}
\usepackage[english]{babel}
\usepackage{csquotes}

\usepackage[nocompress]{cite}

\usepackage{amsmath,amsfonts,amssymb}
\usepackage[standard,amsmath,thmmarks]{ntheorem}
\usepackage{enumitem}

\usepackage{appendix}
\usepackage{marginnote}
\usepackage[hidelinks]{hyperref}

\usepackage{graphicx}
\usepackage[margin=0 cm, font=footnotesize]{caption}

\usepackage[labelformat=simple]{subcaption}

\usepackage{tikz}
\usetikzlibrary{shapes,calc,decorations.pathreplacing,decorations.markings,arrows}
\usetikzlibrary{cd}

\usepackage{mathabx}

\usepackage[lmargin=2.cm,rmargin=2.cm,tmargin=2.5cm,bmargin=2.5cm]{geometry}

\usepackage{todonotes}

\newcommand{\DBF}[1]{\expandafter\newcommand\csname #1\endcsname{{\mathbf{ #1 }}}}

\DBF{T}
\DBF{R}
\DBF{Q}
\DBF{N}
\DBF{Z}
\DBF{C}
\DBF{1}
\DBF{e}
\DBF{I}
\DBF{GL}
\DBF{SL}
\DBF{SO}
\DBF{n}

\newcommand{\DCF}[1]{\expandafter\newcommand\csname cal#1\endcsname{{\mathcal{ #1 }}}}

\DCF{A}
\DCF{C}
\DCF{D}
\DCF{E}
\DCF{H}
\DCF{I}
\DCF{K}
\DCF{M}
\DCF{N}
\DCF{R}
\DCF{B}
\DCF{T}
\DCF{G}
\DCF{O}
\DCF{P}
\DCF{L}
\DCF{S}
\DCF{V}
\DCF{X}
\DCF{Y}
\DCF{F}
\DCF{W}
\DCF{Z}
\DCF{U}
\DCF{Q}

\newcommand{\DSF}[1]{\expandafter\newcommand\csname scr#1\endcsname{{{ #1 }}}}
\DSF{o}
\DSF{w}

\newcommand{\DFF}[1]{\expandafter\newcommand\csname frak#1\endcsname{{\mathfrak{ #1 }}}}

\DFF{L}
\DFF{g}
\DFF{o}
\DFF{so}
\DFF{X}

\newcommand{\DTF}[1]{\expandafter\newcommand\csname #1\endcsname {{\tt #1}}}

\newcommand {\DMO}[1]{\expandafter\DeclareMathOperator\csname #1\endcsname {#1}}

\DMO{diag}
\DMO{supp}
\DMO{id}

\DMO{diam}
\DMO{dom}
\DMO{coker}
\DMO{ind}
\DMO{Lie}
\DMO{sgn}
\DMO{crit}
\DMO{Jac}
\DMO{Hess}
\DMO{grad}
\DMO{rk}
\DMO{codim}
\DMO{SpecFlow}
\DMO{tr}
\DMO{Aut}

\DMO{HC}
\DMO{HTW}

\DMO{Vol}
\DMO{dist}

\DMO{Inv}
\DMO{BInv}

\DMO{div}
\DMO{period}
\DMO{SF}

\mathchardef\hyphen="2D

\newcommand{\setmin}{\smallsetminus}

\newcommand{\iso}{\cong}

\newcommand\rst[2]{ \left. {#1} \right|_{#2} }

\newcommand{\eval}[3]{ \left. {#1} \right|_{#2}^{#3} }

\newcommand{\cl}[1]{ \ensuremath{\overline{#1}} }

\newcommand{\bdy}[1]{\partial {#1} }

\newcommand{\real}{\mathfrak{Re}\,}

\renewcommand{\d}[1]{ \ensuremath{ \operatorname{d}\!{#1} } }

\newcommand{\ddxe}[1]{\ensuremath{ \dfrac{\d{ \hspace{1ex} }}{\d{#1} } }}

\newcommand{\set}[2]{ \left\{ #1 \; : \; #2 \right\} }

\renewcommand{\epsilon}{\varepsilon}
\renewcommand{\phi}{\varphi}

\renewcommand{\leq}{\leqslant}
\renewcommand{\geq}{\geqslant}
\renewcommand{\dots}{\ldots}

\renewcommand{\qed}{\hfill \ensuremath{\Box}}

\makeatletter
\newtheoremstyle{mynonumberplain}%
  {\item[\theorem@headerfont\hskip\labelsep ##1\theorem@separator]}%
  {\item[\theorem@headerfont\hskip \labelsep ##1\
    ##3\theorem@separator]}
\makeatother

\theoremstyle{plain}
\theoremprework{} 
\theorembodyfont{\slshape}
\theorempostwork{}
\theoremseparator{.}
\newcounter{dummy}
\numberwithin{dummy}{section}
\newtheorem{mydef}[dummy]{Definition}
\newtheorem{mythm}[dummy]{Theorem}
\renewtheorem*{mythm*}{Theorem}
\newtheorem{mylemma}[dummy]{Lemma}
\renewtheorem*{mylemma*}{Lemma}
\newtheorem{mycor}[dummy]{Corollary}

\theoremnumbering{Alph}

\theoremnumbering{arabic}
\theoremsymbol{\ensuremath{\qed}}
\newtheorem{myremark}[dummy]{Remark}
\newtheorem{myexample}[dummy]{Example}
\theoremstyle{nonumberplain}
\theoremsymbol{}
\newtheorem{myhyp}{Hypothesis}

\theorembodyfont{\upshape}
\theoremsymbol{\ensuremath{\diamondsuit}}
\theoremseparator{.}
\theoremprework{}

\theoremstyle{mynonumberplain}
\theorembodyfont{\upshape}
\theoremsymbol{\rule{1.2ex}{1.2ex}}
\theoremseparator{.}
\newtheorem{myproof}{Proof}

\theoremstyle{nonumberplain}
\theoremheaderfont{\normalfont\bfseries}
\theoremseparator{.}
\theoremsymbol{}
\theorembodyfont{\slshape}

\setlist[description]{leftmargin=2ex,labelindent=\parindent}

\makeatletter
\newtheoremstyle{mybreak}%
  {\item[\rlap{\vbox{\hbox{\hskip\labelsep \theorem@headerfont
          ##1\ ##2\theorem@separator}\hbox{\strut}}}]}%
  {\item[\rlap{\vbox{\hbox{\hskip\labelsep \theorem@headerfont
          ##1\ ##3,\ continued\theorem@separator}\hbox{\strut}}}]}
\makeatother

\theoremstyle{mybreak}

\makeatletter
\def\namedlabel#1#2{\begingroup
   \def\@currentlabel{\textrm{(#2)}}%
   \label{#1}\endgroup
}
\makeatother

\renewlist{enumerate}{enumerate}{1}
\setlist[enumerate]{label=(\arabic*),}

\begin{document}


\title{Spatial Hamiltonian identities for nonlocally coupled systems\footnotetext{BB is corresponding author, \href{mailto:bente.h.bakker@gmail.com}{\texttt{bente.h.bakker@gmail.com}}}}
\author{Bente Bakker\footnote{Department of Mathematics, VU Universiteit Amsterdam, De Boelelaan 1081, 1081 HV Amsterdam}
\and Arnd Scheel\footnote{School of Mathematics, University of Minnesota, 206 Church St. SE, Minneapolis, 55455}
}
\date{\today}
\maketitle

\begin{abstract}
We consider a broad class of systems of nonlinear integro-differential equations posed on the real line that arise as Euler-Lagrange equations to energies involving nonlinear nonlocal interactions. Although these equations are not readily cast as dynamical systems, we develop a calculus that yields a natural Hamiltonian formalism. In particular, we formulate Noether's theorem in this context, identify a degenerate symplectic structure, and derive Hamiltonian differential equations on finite-dimensional center manifolds when those exist. Our formalism yields new natural conserved quantities. For Euler-Lagrange equations arising as traveling-wave equations in gradient flows, we identify Lyapunov functions. We provide several applications to pattern-forming systems including neural field and phase separation problems. 
\end{abstract}

\paragraph{\footnotesize 2010 Mathematics Subject Classification.}
{\footnotesize
Primary: 35S30, 45G15; \hspace{1 pt} Secondary: 35C07, 37K05, 37L10, 37L45.
}

\paragraph{\footnotesize Acknowledgements.}
{\footnotesize
BB is partially supported by NWO VICI grant 639.033.109 and NWO TOP grant 613.001.351.

AS gratefully acknowledges support through NSF grant DMS--1311740.
}


\binoppenalty=\maxdimen
\relpenalty=\maxdimen
\sloppy


\section{Introduction}

Differential equations of gradient form $u_t = - \nabla \calE(u)$ or of Hamiltonian form $u_t = J \nabla \calE(u)$
arise throughout mathematical modeling from maximal energy dissipation first principles, or from least action principles, respectively.
Mathematically, the energy $\calE$ is simply a function defined over a finite or infinite-dimensional space. The associated gradient and Hamiltonian flows are  differential equations with equilibria given by the critical points of the energy. Energies over infinite-dimensional spaces are usually defined on function spaces through integrals of nonlinear functions of state variables and their derivatives. Critical points then solve Euler-Lagrange equations, commonly of elliptic type. Dependence of the energy on derivatives of the state variables encodes \emph{local} interactions, often derived in various types of continuum limits. We are interested here in cases where this continuum limit retains \emph{nonlocal} interaction terms. More specifically, we are interested in the somewhat specific class of energies $\calE$ that contain nonlocal interaction term of the form
\[
\iint S(u(x))^T K(x-y) S(u(y)) \d y \d x, \qquad x,y\in\R, \quad u(\cdot) \in\R^d,  
\]
modeling phenomena such as long-range interactions of agents in social interaction, between particles through nonlocal force fields, or of neurons, labeled in a feature space $x$  through synaptic connections. In all those cases, the convolution structure embodies the modeling assumption of translational invariance of physical space. 
We are concerned here with critical points of the energy $\calE$, possibly up to actions of the underlying translation symmetry group. Such critical points, interpreted as solutions in physical space, are referred to as solitary waves, excitation pulses, waves trains, or traveling fronts, depending on the modeling context.
We refer to \cite{bates1997traveling, bates1999integrodifferential, bates2006some, du2012analysis, du2013analysis, silling2000reformulation} for examples from material sciences,
\cite{ehrnstrom2009traveling, ehrnstrom2012existence, naumkin1994nonlinear, whitham2011linear} for examples from fluid dynamics,
and \cite{ermentrout1993existence, pinto2001spatially, wilson1972excitatory, wilson1973mathematical} for examples from neurobiology.

For local equations, dynamical systems methods have provided powerful tools to the study of such coherent structures. One therefore casts the Euler-Lagrange equations, which are simply systems of higher-order ordinary differential equations,  as dynamical systems in the spatial variable $x$, which is then considered as a time-like variable, trying to describe the set of bounded solutions, including periodic, heteroclinic, and homoclinic solutions. The dynamical systems tools can then be thought of as \emph{pointwise} methods, exploiting the geometry of the ``spatial phase space'' through tools including  phase plane analysis with the Poincar\'e-Bendixson theorem, center manifold reduction and normal form techniques, bifurcation and geometric singular perturbation theory, or topological index theory.

For nonlocal equations, much of this pointwise technology is not immediately available. As the main obstacle, nonlocal interactions in space typically generate forward- and backward delay terms in the time-like spatial variable and it is often not clear how a pointwise formulation as a dynamical system in a phase space can be recovered. Much of the previous work has therefore focused on variational or perturbative techniques, although even those may sometimes benefit from pointwise estimates. More recently, certain dynamical systems techniques  have been made available for the study of nonlocal equations. A first approach casts systems with finite-range interactions as ill-posed dynamical systems on an infinite-dimensional phase space, much like elliptic equations posed on a cylinder, and uses dynamical systems techniques, in particular invariant manifolds, to study the dynamics; see for instance \cite{hss,mpvl,hvl,hupkes}. A different avenue evokes dynamical systems techniques without ever setting up a phase space for a pointwise description, but rather adapting techniques from dynamical systems by reducing to basic functional analytic aspects; see \cite{faye2015existence}  for  geometric singular perturbation techniques, \cite{faye2016center}  for center manifold reductions,
and \cite{scheel2017bifurcation} for bifurcation methods. 

The present work can be viewed as continuing this latter approach, providing a framework for exploiting Hamiltonian and symplectic concepts. Indeed, for pointwise functionals, the Euler-Lagrange equations inherit a Hamiltonian structure, viewing the integral of the energy density as an action functional, and the energy density as the Lagrangian. Using Legendre transform in the case of dependence on gradients, only, one readily finds the associated Hamiltonian for the Euler-Lagrange equation. Such a procedure is slightly more involved for functionals containing higher-order derivatives. We show how to formulate the Euler-Lagrange equations in a symplectic framework, identifying symplectic structure, Hamiltonian and more general conserved quantities starting from Noether's theorem and a non-local integration-by-parts calculation. In the case of traveling-wave equations, conserved quantities turn into Lyapunov functions. 

Hamiltonian formalisms for infinite-dimensional systems are of course not new. They had been exploited also for ill-posed, elliptic equations; see for instance \cite{kirch,mielke2006hamiltonian,gui2008hamiltonian}. For nonlocal equations based on the fractional Laplacian, an extension to an additional spatial dimension using the Dirichlet-to-Neumann operator on the half space allows one to cast the nonlocal equation as a local equation, again obtaining conserved quantities through a local Hamiltonian formalism \cite{cabre,sire}. We are not aware of cases where such identities have been derived in the presence of ``truly'' nonlocal interactions. 

Besides providing a systematic and elegant framework to formulate equations, a Hamiltonian formalism provides useful tools for the analysis. We therefore provide a somewhat lengthy list of applications that illustrate how the formalism here can be put to work. 
In these examples $K : \R \to \calM(\R^d)$ will always be a curve of symmetric matrix convolution kernels,
that is, $K(r) = K(-r)$ and $K(r) = K(r)^T$,
with sufficiently rapid decay as $r \to \pm \infty$.
\begin{myexample}[Nonlocal Allen-Cahn]
\label{ex:Allen-Cahn}
Let $F : \R^d \to \R$ be a differentiable map, and consider the energy
\[
\calE(u) := \frac 1 4 \iint \big( u(x)-u(y) \big)^T K(x-y) \big( u(x) - u(y) \big) \d y \d x + \int F(u(x)) \d x.
\]
Then the formal gradient flow $u_t = - \nabla_{L^2} \calE(u)$ is the nonlocal analogue of the Allen-Cahn equation, as introduced in \cite{bates1997traveling}
as a model for phase transition,
\[
u_t = - \kappa_0 u + \int K(x-y) u(t,y) \d y - \nabla F(u),
\]
where $\kappa_0 = \int K(x) \d x$.
Typically, the interaction kernel $K$ in these models is algebraically localized and positive, and the potential $F$ exhibits multiple local minima. 
With the traveling wave ansatz $u(t,x) = u(\xi)$, $\xi = x - ct$, we obtain the integro-differential equation
\begin{equation}
  \label{eq:Allen-Cahn_TW}
  -c u' = - \kappa_0 u + K*u - \nabla F(u).
\end{equation}
\end{myexample}

\begin{myexample}[Neural field equation]
\label{ex:NFE}
Neural field equation, introduced in \cite{wilson1972excitatory, wilson1973mathematical},
\begin{equation}
  \label{eq:NFE}
  u_t = -u + \int K(x-y) S(u(t,y)) \d y,
\end{equation}
model interactions of neurons through nonlocal synaptic connections. Typically, one thinks of $x$ as a feature space, grouping neurons by function and/or location, and $S : \R^d \to \R^d$ as a smooth sigmoidal input function, triggering firing of a neuron above a threshold value. The sign of the interaction kernel $K$ encodes inhibitory versus excitatory connections. Sign changes in $K$ can lead to complex stationary patterns and have been proposed as a template for short-term memory functions.  

In the specific case that $DS(u)$ is positive definite for all $u \in \R^d$, we assume that there exists a map $E:\R^d \to \R$ so that $\nabla E(u) = DS(u)u$. We then consider the energy
\[
\calE(u) := \int E(u(x)) \d x - \frac 1 2 \iint S(u(x))^T K(x-y) S(u(y)) \d y \d x,
\]
and write equation \eqref{eq:NFE} as 
\[
u_t = - \big( DS(u)^T \big)^{-1} \nabla_{L^2} \calE(u).
\]
Hence, the neural field equation can formally be written as the gradient flow of $\calE$ on a Hilbert manifold $\calM$ modeled over $L^2$,
where the Riemannian metric is given by
\[
g_u(v,w) = \langle v , DS(u) w \rangle_{L^2}, \qquad \text{where} \quad v, w \in T_u\calM.
\]
With the traveling wave ansatz $u(t,x) = u(\xi)$, $\xi = x - ct$ in \eqref{eq:NFE}, we obtain
\[
-c u' = - u + K*S(u).
\]
\end{myexample}

In both of these examples the original equation is a gradient flow. By analogy to examples from local partial differential equations, one expects in such a situation that  the equations governing traveling waves have a gradient-like structure.
In a slight alteration, the next example starts out with a Hamiltonian system.
The traveling wave equation now inherits a variational structure, and  one expects to find a Hamiltonian system from the Euler-Lagrange equations.

\begin{myexample}[Whitham type equation]
\label{ex:Whitham}
Consider the energy
\[
\calE(u) := \frac \alpha 3 \int u^3 \d x + \frac 1 2 \iint u(x) K(x-y) u(y) \d y \d x,
\]
where $\alpha \in \R$ is a parameter and $u : \R \to \R$.
The Hamiltonian equation $u_t = J \nabla_{L^2} \calE(u)$ with $J = - \partial_x$ is Whitham's equation, introduced in \cite{whitham2011linear} as a model for shallow water waves with weak dispersion,
\begin{equation}
  \label{eq:Whitham}
  u_t + 2 \alpha u u_x + \int K(x-y) u_y(t,y) \d y = 0,
\end{equation}
where
\[
K(x) = \int_\R e^{i x \xi} \sqrt{ \frac{ \tanh(\xi) }{ \xi } } \d \xi.
\]
This model incorporates a generic second-order nonlinearity as found in the KdV approximation, 
but replaces the dispersive terms in the KdV equation by a convolution operator which models the full frequency dispersion found in 
linear surface water waves in finite depth.
Making the traveling wave ansatz $u(t,x) = u(\xi)$, $\xi = x - c t$, equation \eqref{eq:Whitham} becomes
\[
\bigg( \frac \alpha 2 u^2 - c u + K * u \bigg)' = 0,
\]
hence
\[
\alpha u^2 - c u + K * u = \mu,
\]
for some $\mu \in \R$.
\end{myexample}

A shared property of these first three examples is invariance under the shift map $u \mapsto u(\cdot + \tau)$, a symmetry which we will see gives rise to a Hamiltonians and Lyapunov functions. This translation symmetry is broken in the next example, which however does possess a rotation symmetry. In this case, the Hamiltonian formalism gives rise to a conserved quantity analogous to an angular momentum.

\begin{myexample}[Nonlocal NLS]
\label{ex:triggered_GL}
\label{ex:NLS}
Let  $S(u) = f( |u|^2 )$, $f : \R \to \R$ differentiable, $\lambda \in \R$, and consider the energy
\[
\calE(u) := - \frac 1 2 \int |u_x|^2 \d x - \frac \lambda 2 \int |u|^2 + \frac 1 2 \iint S(u(x)) K(x-y) S(u(y)) \d y \d x,
\]
and extend $\nabla_{L^2} \calE(u)$ linearly to complex valued functions $u : \R \to \C$.
The Hamiltonian flow $u_t = i \nabla_{L^2} \calE(u)$ is the nonlocal nonlinear Schr\"odinger (nonlocal NLS) equation,
\begin{equation}
  \label{eq:NLS}
  - i u_t = u_{xx} - \lambda u + DS(u) K*S(u).
\end{equation}
We study here the effect of a (propagating) parameter jump type inhomogeneity $\lambda=\lambda(x-ct)$ which,  after rescaling, can be assumed to satisfy
\[
\lambda(\xi) =
\begin{cases}
  0 \quad &\text{if } \xi \leq -\ell, \\
  1 \quad &\text{if } \xi \geq \ell,
\end{cases}
\]
for some value of $\ell  > 0$.
Substituting the ansatz $u(t,x) = e^{i c x / 2} A(x-c t)$ in \eqref{eq:NLS}, we obtain
\[ 
A'' + \left( \lambda(\xi) - \frac{c^2}{4} \right) A + DS(A) K* S(A) = 0.
\]
Note that, assuming $S$ is at least quadratic, the medium supports a band of small-amplitude plane wave solutions only for $\xi\geq \ell$, when $0<|c|<2$, but not for $\xi\leq -\ell$.  We will study here how the parameter jump affects this band near $\xi=+\infty$. We refer to \cite{gui2008hamiltonian,lloyd,monteiro,weinburd} for more (elaborate) examples of wavenumber (or angle) selection at interfaces.
\end{myexample}

Note that in these examples we dealt with the energy functionals, gradient flows and Hamiltonian flows 
only formally in that integrals may not converge on typical solutions of interest.
%
%
%

\paragraph{Outline of main results.}
We describe our main results informally and outline the remainder of the paper. 
We are concerned with variational problems of the form
\[
\int L\bigg(x,u(x),u'(x),K*S(u)(x)\bigg) \d x \quad \leadsto \quad \text{minimize},\]
and the associated Euler-Lagrange equations. We are interested in conserved quantities, that is, in nontrivial maps $\calC : C^1(\R,\R^d) \to C^0(\R,\R^d)$ that map solutions of the Euler-Lagrange equations to constant functions. We therefore establish a nonlocal analogue of Noether's theorem, allowing us to derive conserved quantities from symmetries of the Euler-Lagrange equation. Key ingredient to the proof is a nonlocal  ``integration by parts'' formula,
\[
\int_a^b A u \cdot v \d x = \eval{\calB(u(\cdot + \tau))}{\tau=a}{b} + \int_a^b u \cdot A^t v \d x,
\]
where $A$ is a nonlocal operator, $A^t$ denotes its formal adjoint, and $\calB$ should be interpreted as ``boundary data''. Specific assumptions and statements, as well as limitations in terms of regularity are contained in Section \ref{sec:Noether}.

In Section \ref{sec:Hamiltonian}, we consider translation invariant  Euler-Lagrange equations
\begin{equation}
  \label{eq:prototype_Hamiltonian}
  - \partial_x \nabla_v(u,u') + \nabla_u E(u,u') + DS(u)^T K*S(u) = 0.
\end{equation}
Through our Noether theorem, the shift invariance leads to a conserved quantity which is itself equivariant with respect to the shift action.
This naturally leads to a Hamiltonian, that is, a map $\calH : C^1(\R,\R^d) \to \R$ such that
\begin{equation}
  \label{eq:prototype_conserved_quantity}
  \ddxe{\tau} \calH\big(u(\cdot + \tau)\big) = 0 \qquad \text{for any solution $u$ of \eqref{eq:prototype_Hamiltonian},}
\end{equation}
and a (degenerate) symplectic form $\omega$ which
allows us to conclude, at least formally,  that the shift action on the set of solutions of the Euler-Lagrange equation is the Hamiltonian flow of $\calH$.

\begin{myremark}
  We reiterate that we are considering Hamiltonian and symplectic structures in the \emph{spatial} variable.
To appreciate this, take, for example, the nonlocal wave equation
\begin{equation}
  \label{eq:Nwave}
  u_{tt} = - u + K*u + f(u), \qquad \big( u(t,x), \; u_t(t,x) \big) \to (0,0) \quad \text{as} \quad x \to \pm\infty.
\end{equation}
This equation has a  \emph{temporal shift} symmetry, and the dynamics are Hamiltonian with respect to the canonical symplectic structure.
Since \eqref{eq:Nwave} also has a \emph{spatial shift} symmetry, the momentum
\[
p(u) = \int_\R u_t u_x \d x,
\]
is conserved by the \emph{temporal dynamics}.
This conserved momentum $p$ is not the conserved quantity described in \eqref{eq:prototype_conserved_quantity}, 
as we are concerned with quantities which are conserved under \emph{spatial shifts} of the stationary solutions.
The interplay between spatial and temporal Hamiltonian structures has been further pursued, at least for local equations, 
using multisymplectic structures, see \cite{bridges1997multi, bridges2001multi, bridges2004nonlinear}.
\end{myremark}

To give a specific example here,  equation \eqref{eq:Allen-Cahn_TW} with $c = 0$
is Hamiltonian with 
\[
    \calH(u) = \rst{ \bigg( F(u) + \frac 1 2 u \cdot \kappa_0 u - \frac 1 2 u \cdot K * u \bigg) }{ x=0 } 
+ \frac 1 2 \int_{x<0} \int_{y>0} u(x) \cdot K(x-y) u'(y) - u'(x) \cdot K(x-y) u(y) \d y \d x,
\]
whilst the symplectic form reads
\[
\omega(v,w) = \int_{x < 0} \int_{y > 0} v(x) \cdot K(x-y) w(y) - w(x) \cdot K(x-y) v(y) \d y \d x.
\]
Inspecting this example, one quickly notes that $\omega$ is in general degenerate:
if the kernel has support $\supp K \subset [-1,1]$ and $v$ has support $\supp v \subset \R \setmin [-1,1]$, then $\omega( v , w ) = 0$, for any $w$. We are however able to show that $\omega$ is locally non-degenerate on a center manifold, giving a (non-canonical) symplectic structure and reduced Hamiltonian vector fields and thereby  extending the results from \cite{faye2016center}.

Section \ref{sec:grad_like} considers the effect of (spatially) dissipative terms, such as equations of the form
\begin{equation}
  \label{eq:prototype_dissipative}
  - \partial_x \nabla_v(u,u')  + \nabla_u E(u,u') + DS(u)^T K*S(u) = \Gamma(u) u',
\end{equation}
where $\Gamma : \R^d \to \calM(\R^d)$ takes it values in the cone of positive definite matrices.
Shift invariance now  leads to a Lyapunov function, $\calL : C^1(\R,\R^d) \to \R$,
such that
\[
\ddxe{\tau} \calL\big(u(\cdot+\tau)\big) \leq 0 \qquad \text{for any solution $u$ of \eqref{eq:prototype_dissipative},}
\]
with equality only for constant $u$. A specific example is  \eqref{eq:Allen-Cahn_TW} , with $c \neq 0$, and
\[
    \calL(u) = \rst{ \bigg( F(u) + \frac 1 2 u \cdot \kappa_0 u - \frac 1 2 u \cdot K * u \bigg) }{ x=0 } 
+ \frac 1 2 \int_{x<0} \int_{y>0} u(x) \cdot K(x-y) u'(y) - u'(x) \cdot K(x-y) u(y) \d y \d x.
\]
As a consequence, all bounded solutions are either constant or heteroclinic with respect to the shift action.

Finally, Section \ref{sec:applications} discusses several applications, elaborating on the four examples presented here in the introduction. We prove a local bifurcation results for small-amplitude periodic and homoclinic solutions of  example \ref{ex:Allen-Cahn} using Hamiltonian center-manifold reduction; we establish the existence of traveling fronts for example \ref{ex:NFE} exploiting a Lyapunov function and Conley index theory, extending results in  \cite{ermentrout1993existence} where positivity of $K$ was required; and we establish local bifurcation of periodic and solitary waves for  example \ref{ex:Whitham}, providing alternative geometric proofs for results in \cite{ehrnstrom2009traveling, ehrnstrom2012existence}. Finally, we establish selection of wavenumbers through parameter jumps in example \ref{ex:triggered_GL}.




\section{Symmetries and conserved quantities}
\label{sec:Noether}
We formulate a setup and state our nonlocal version of Noether's theorem. 
\subsection{Variational setup}
  
\paragraph{Nonlocal variational problems.}
We will throughout rely on the following assumption on the convolution kernel.  $K : \R \to \calM(\R^d)$.
\begin{myhyp}[K]
\namedlabel{hyp:K}{K} The map $K$ is symmetric, meaning $K(r) = K(-r)$ and $K(r) = K(r)^T$.
We assume $K \in L^1(\R,\calM(\R^d))$, that is, $K_{ij} \in L^1(\R,\R)$ for all $1 \leq i,j \leq d$.
Moreover,
\[
\int_\R \left( 1 + |r| \right) \| K(r) \| \d r < \infty.
\]
\end{myhyp}
Now let $S \in \C^2( \R^d ,\R^d)$ and define the nonlinear operator
\[
\calN(u)(x) := K*S(u)(x) = \int_\R K(x-y) S(u(y)) \d y.
\]
Then $\calN$ is (formally) differentiable, with
\[
D \calN(u)[v](x) = K*(DS(u) v)(x) = \int_\R K(x-y) DS(u(y)) v(y) \d y,
\]
and formal $L^2$-adjoint 
\[
D \calN(u)^t[w](x) := DS(u(x))^T K*w(x) = DS(u(x))^T \int_\R K(x-y) w(y) \d y.
\]
Let the Lagrangian  $L \in C^2(\R \times \R^d \times \R^d \times \R^d,\R)$ be given as a function $L=L(x,u,v,n)$;
from Section \ref{sec:Hamiltonian} onward we will impose further restrictions on $L$.
For $a,b \in \R$, the ``truncated action'' ${\calA_a^b : C^2_b(\R,\R^d) \to \R}$ is defined through
  \[
  \calA_a^b(u) := \int_a^b L\big(x,u(x),u'(x),K*S(u)\big) \d x.
  \]
Here $C^k_b(\R,\R^d)$ denotes the space of $C^k$ functions whose first $k$ derivatives are uniformly bounded.
Note that for finite $a$, $b$ the integral is convergent. We shall also write $\calA$ to denote the map ${(u,a,b) \mapsto \calA_a^b(u)}$.

Now consider the nonlocal equation
\begin{equation}
  \label{eq:EL}
  \nabla_u L(x,u,u',K*S(u)) - \partial_x \nabla_v L(x,u,u',K*S(u)) +DS(u)^T K* \nabla_n L(x,u,u',K*S(u)) = 0.
\end{equation}
  Note that, formally, equation \eqref{eq:EL}
  is the Euler-Lagrange equation of the action functional $\calA_{-\infty}^\infty$.
  We will not attempt to make this statement rigorous, since the integral will typically not converge.

\begin{mydef}[Regularity of  solutions]
\label{def:solution_regularity}
  We will consider solutions of \eqref{eq:EL} which are of class $C^R$, where $R = 1$ if $\nabla_v L = 0$, and $R = 2$ otherwise.
\end{mydef}

  \begin{myremark}
    If $\nabla_v L \neq 0$, taking $R = 2$ corresponds to considering strong solutions of \eqref{eq:EL}.
    On the other hand, if $\nabla_v L = 0$, \eqref{eq:EL} is well-defined on functions of class $L^\infty$.
    However, the proof of Theorem \ref{thm:Noether} fails for $u \in L^\infty$ without further regularity assumptions.
    In fact, the conclusion of the Noether theorem is in general not valid for solutions of \eqref{eq:EL} lacking regularity; see Remark \ref{r:dis}.
    We therefore restrict our attention to solutions of class $C^1$.
    We will give examples of bootstrap arguments that establish the required regularity in some cases in Sections \ref{sec:Hamiltonian} and \ref{sec:grad_like}.
  \end{myremark}
  \begin{myremark}
    When $R = 1$, it suffices to consider $L$ and $S$ of class $C^1$.
    We will then extend the domain of $\calA_a^b$ to $C^1_b(\R,\R^d)$.
  \end{myremark}

\paragraph{Symmetry.}
Let $G$ be a Lie subgroup of $\GL(d) \times \R$, acting on functions $u : \R \to \R^d$ via the canonical representation
\[
g \bullet u := \phi_\tau \big( A u \big), \qquad \text{for} \quad g = (A,\tau) \in G.
\]
Here $\phi_\tau$ is defined by
\[
\phi_\tau u(x) := u(x + \tau), \qquad \text{for all} \quad x \in \R.
\]
By slight abuse of notation, we will use the same notation for the action of $G$ on various function spaces.
We also let $G$ act on $\R$ by
\[
g \bullet x := x - \tau, \qquad \text{where} \quad g = (A,\tau) \in G, \quad x \in \R.
\]

  Let $\frakg \iso T_e G$ be the Lie algebra of $G$.
  Formally, define the induced action on functions $u : \R \to \R^d$ by
  \[
  u_\xi := \rst{ \ddxe{\tau} \exp_e( \tau \xi ) \bullet u }{\tau=0}, \qquad \text{where} \quad \xi \in \frakg.
  \]
  Note that, when $G$ contains the translations $\{ \I_d \} \times \R$, we have $u_\xi \in C^{R-1}_b(\R,\R^d)$ whenever $u \in C^R_b(\R,\R^d)$.
  Similarly, we define the induced action on $\R$ by
  \[
  \1_\xi := \rst{ \ddxe{\tau} \exp_e( \tau \xi ) \bullet 1 }{\tau=0}, \qquad \text{where} \quad \xi \in \frakg.
  \]
  
  \begin{mydef}
  \label{def:symmetry}
  A group $G$ is called the global symmetry group of \eqref{eq:EL} if it is the maximal Lie subgroup of $\GL(d) \times \R$ for which
  \[
  \calA\big( g \bullet (u,a,b) \big) = \calA\big( u,a,b \big), \qquad \text{for all} \quad (u,a,b) \in C^R_b(\R,\R^d) \times \R^2, \quad g \in G.
  \]
  Similarly, the Lie algebra $\frakg$ is called the local symmetry algebra of \eqref{eq:EL} if it is maximal and
  \[
  D \calA( u , a , b )(u_\xi,\1_\xi,\1_\xi) = 0, \qquad \text{for all} \quad (u,a,b) \in C^R_b(\R,\R^d) \times \R^2, \quad \xi\in\frakg.
  \]
  \end{mydef}

\subsection{Noether's theorem for nonlocal equations}

Noether's theorem in the classical setting relates symmetries and conserved quantities; see for example \cite{marsden2013introduction}. We extend this result here to our class of nonlocal equations \eqref{eq:EL}.

\begin{figure}
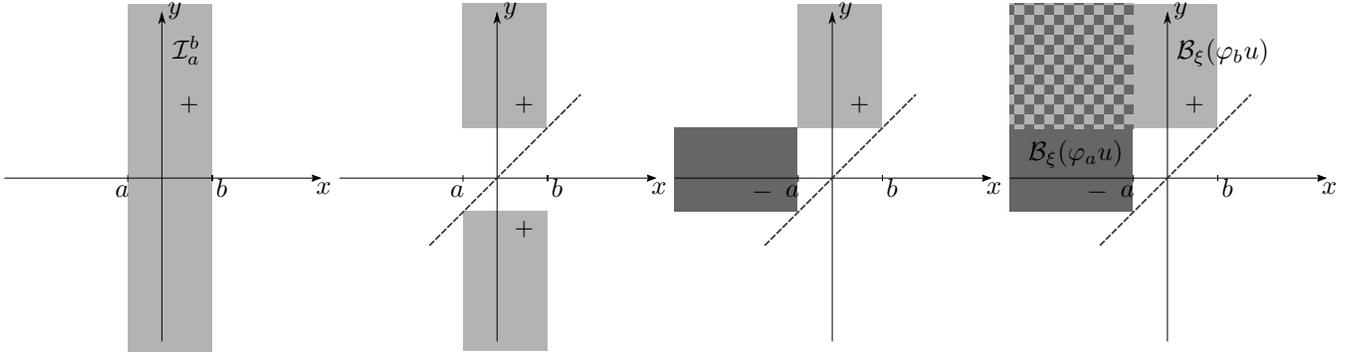

  \centering
  \begin{subfigure}[c]{0.24\textwidth}
    \def\svgwidth{\textwidth}
    \input{fig/int1.tex}
  \end{subfigure}
  \hfill 
  \begin{subfigure}[c]{0.24\textwidth}
    \def\svgwidth{\textwidth}
    \input{fig/int2.tex}
  \end{subfigure}
\hfill
  \begin{subfigure}[c]{0.24\textwidth}
    \def\svgwidth{\textwidth}
    \input{fig/int3.tex}
  \end{subfigure}
  \hfill
  \begin{subfigure}[c]{0.24\textwidth}
    \def\svgwidth{\textwidth}
    \input{fig/int4.tex}
  \end{subfigure}
  \caption{Transformation of $\calI_a^b := \iint_{\Omega_0} s(x,y) \d V$ into boundary terms $\eval{ \calB_\xi(\phi_\tau u) }{\tau=a}{b}$. From left to right, we show the original domain of integration, then the square $(a,b)^2$ can be removed by skew symmetry of the integrand, then, using the skew symmetry, the domain of integration is transformed. On the right, the resulting difference of integrals coincides with the difference of the boundary terms, adding and subtracting the checkered region.}
  \label{fig:computing_boundary_terms}
\end{figure}

\begin{mythm}[Noether]
\label{thm:Noether}
  Let $\frakg$ be the local symmetry algebra of \eqref{eq:EL}.
  For any given $\xi \in \frakg$, we define the function ${\calC_\xi :C^R_b(\R,\R^d) \to C^0_b(\R,\R)}$ by
  \[
  \calC_\xi(u)(x) := L\big(x,u(x),u'(x),K*S(u)(x)\big) \1_\xi + \nabla_v L(x,u,u',K*S(u)) \cdot u_\xi + \calB_\xi(\phi_x u).
  \]
  Here
  \begin{equation}\label{e:Q}
  \calB_\xi(u) := \iint_{Q} \sigma(x,y) - \sigma(y,x) \d x \d y, \qquad \text{with} \quad  Q := (-\infty,0) \times (0,\infty)
  \end{equation}
  and
  \[
    \sigma(x,y) := \rst{\nabla_n L\big(x,u,u',K*S(u)\big)}{x} \cdot  K(x-y) \rst{ \big( DS(u) u_\xi \big) }{ y }.
  \]
  Then $x \mapsto \calC_\xi(u)(x)$ is constant,
  for any $u \in C^R_b(\R,\R^d)$ which solves \eqref{eq:EL}.
\end{mythm}
\begin{myproof}
  Note that the integrand of $\calA_a^b(\exp_e(\tau \xi) \bullet u)$ is $C^1$ with respect to $\tau$,
  hence we can differentiate under the integral sign.
  We find
  \begin{equation}
    \label{eq:step1}
   \begin{split}
    0 &= D \calA(u,a,b)(u_\xi,a_\xi,b_\xi) = \rst{ \ddxe{\tau} \calA( \exp_e(\tau \xi) \bullet (u,a,b) )}{\tau=0} \\
        &= \eval{ L\big(x,u(x),u'(x),K*S(u)(x)\big) \1_\xi }{x=a}{b} \\
    &\quad + \int_a^b \nabla_u L(x,u,u',K*S(u)) \cdot u_\xi + \nabla_v L(x,u,u',K*S(u)) \cdot u_\xi' \d x \\
    &\quad + \int_a^b \nabla_n L\big(x,u,u',K*S(u)\big) \cdot K*DS(u) u_\xi \d x.
  \end{split} 
  \end{equation}
  Upon integration by parts, we find
  \begin{multline}
    \label{eq:step2}
    \int_a^b \nabla_v L(x,u,u',K*S(u)) \cdot u_\xi' \d x = \eval{ \nabla_v L(x,u,u',K*S(u)) \cdot u_\xi }{x=a}{b} 
    - \int_a^b \partial_x \nabla_v L(x,u,u',K*S(u)) \cdot u_\xi \d x.
  \end{multline}
  Inserting this back into \eqref{eq:step1}, and using that $u$ solves \eqref{eq:EL}, we find
  \begin{multline}
    \label{eq:noether_step}
    0 = \eval{ L\big(x,u(x),u'(x),K*S(u)(x)\big) \1_\xi }{x=a}{b} + \eval{ \nabla_v L(x,u,u',K*S(u)) \cdot u_\xi }{x=a}{b} \\
    + \int_a^b \nabla_n L\big(x,u,u',K*S(u)\big) \cdot K * DS(u) u_\xi
    - DS(u)^T K*\nabla_n L\big(x,u,u',K*S(u)\big) \cdot u_\xi \d x.
  \end{multline}
  
  We will now express the last integral as a boundary term, that is, we need to find a map $\calB_\xi$ such that
  \begin{align*}
      \int_a^b \nabla_n L\big(x,u,u',K*S(u)\big) \cdot K * DS(u) u_\xi
    - DS(u)^T K*\nabla_n L\big(x,u,u',K*S(u)\big) \cdot u_\xi \d x \\
    = \calB_\xi(\phi_b u) - \calB_\xi(\phi_a u).
  \end{align*}
  This can be interpreted as a nonlocal analogue of Green's formula, similar to the nonlocal vector calculus developed in \cite{gunzburger2010nonlocal}.
  We note that the cited work cannot be applied directly, because the integrand
  \[
  \nabla_n L\big(x,u,u',K*S(u)\big) \cdot K * DS(u) u_\xi
    - DS(u)^T K*\nabla_n L\big(x,u,u',K*S(u)\big) \cdot u_\xi 
  \]
  is not integrable over $\R$.
  We note that
  \begin{multline*}
\int_a^b \nabla_n L\big(x,u,u',K*S(u)\big) \cdot K*DS(u) u_\xi
    - DS(u)^T K* \nabla_n L\big(x,u,u',K*S(u)\big) \cdot u_\xi \d x \\     
= \int_a^b \int_\R \sigma(x,y) - \sigma(y,x) \d y \d x.
  \end{multline*}
  The geometric rationale for the transformation of the integral into boundary terms is depicted in Figure \ref{fig:computing_boundary_terms}.
  Let $s(x,y) := \sigma(x,y) - \sigma(y,x)$.
  We now decompose the iterated integral into three integrals over the domains
  $\Omega_a = (a,b) \times (-\infty,a)$, $\Omega_0 = (a,b) \times (a,b)$, and $\Omega_b = (a,b) \times (b,\infty)$.
  Let $\Psi : \R^2 \to \R^2$ be the reflection around the diagonal: $\Psi(x,y) = (y,x)$.
  Then $\Psi^{-1}(\Omega_0) = \Omega_0$ and $s \circ \Psi = - s$, hence
  \[
  \iint_{\Omega_0} \sigma(x,y) - \sigma(y,x) \d x \d y = \iint_{\Psi^{-1}(\Omega_0)} s \circ \Psi \d x \d y = - \iint_{\Omega_0} s \d x \d y = 0.
  \]
  Moreover, if we let $Q_\tau := (-\infty,\tau) \times (\tau,\infty)$, then
  \[
  \Psi^{-1}(\Omega_a) = Q_a \setmin (Q_a \cap Q_b), \qquad \Omega_b = Q_b \setmin (Q_a \cap Q_b).
  \]
  The conditions imposed on $K$ in Hypothesis \ref{hyp:K} ensure that $s$ is integrable over $Q_\tau$, for each $\tau \in \R$.
  Indeed, reparameterizing $Q_0$ by
  \[
  \begin{array}{l l}
    x = ( r + s ) / 2, \qquad & r < 0, \\
    y = ( s - r ) / 2, \qquad & r < s < -r,
  \end{array}
  \]
  we find
  \begin{equation}
    \label{eq:K_moment_Q_integrable}
    \begin{split}
      \iint_{Q_\tau} \|K(x-y)\| \d x \d y &= \iint_{Q_0} \|K(x-y)\| \d x \d y = \frac 1 2 \int_{r<0} \int_{|s|<|r|} \| K(r) \| \d s \d r \\
      &= \int_{r < 0} |r| \| K(r) \| \d r = \frac 1 2 \int_\R |r| \|K(r)\| \d r,
    \end{split}
  \end{equation}
  which is finite by assumption.
  
  We now have
  \[
  \iint_{\Omega_a} s \d x \d y = \iint_{Q_a} s \circ \Psi \d x \d y - \iint_{Q_a \cap Q_b} s \circ \Psi \d x \d y = - \iint_{Q_a} s \d x \d y + \iint_{Q_a \cap Q_b} s \d x \d y,
  \]
  and 
  \[
  \iint_{\Omega_b} s \d x \d y = \iint_{Q_b} s \d x \d y - \iint_{Q_a \cap Q_b} s \d x \d y.
  \]
  Therefore,
  \[
  \iint_{(a,b)\times\R} s \d x \d y = \iint_{\Omega_a} s \d x \d y + \iint_{\Omega_b} s \d x \d y = \iint_{Q_b} s \d x \d y - \iint_{Q_a} s \d x \d y.
  \]
  Then note that
  \[
  \iint_{Q_\tau} s \d x \d y = \iint_{Q_0} s(x+\tau,y+\tau) \d x \d y = \calB_\xi( \phi_\tau u ).
  \]
  The last identity follows because the action of $\xi$ commutes with the $\R$-action of $\phi_\tau$.
  After collecting the various terms, the claimed result follows.
\end{myproof}

\begin{myremark}\label{r:dis}  We note here that the Noether theorem is in general not valid for discontinuous solutions.
  As an  example, consider  the stationary nonlocal Allen-Cahn equation
  \begin{equation}
    \label{eq:pinning}
    0 = d ( - u + K* u ) + f_a(u),
  \end{equation}
  where $K$ is a positive, smooth kernel with $\int K(x) \d x = 1$, $d > 0$ is a coupling constant, and $f_a(u) = u ( 1 - u )  ( u - a )$.
  The results from \cite{anderson2016pinning, bates1997traveling} show that for $0 < d < 1/4$, there exist open intervals $\left(a_-(d),a_+(d)\right) \subset (0,1)$
  such that the following holds.
  For any $a \in \left(a_-(d),a_+(d)\right)$, equation \eqref{eq:pinning} has a solution $u$ with precisely one jump discontinuity, and such that $\lim_{x\to -\infty} u(x) = 1$, $\lim_{x\to\infty} u(x) = 0$.
  The conserved quantity corresponding to the translational symmetry in \eqref{eq:pinning} is given by $\calC_\xi(u)(x) = \calH(\phi_x u)$, where
  \[
      \calH(u) = \rst{ \bigg( \frac d 2 u^2 - F_a(u) - \frac d 2 u K * u \bigg) }{ x=0 } 
    + \frac d 2 \iint_Q K(x-y) \bigg( u(x)  u'(y) - u'(x) u(y) \bigg) \d x \d y,
  \]
  with $F_a(u) = \int_0^u f_a(s) \d s$.
  But $\calH(0) = - F_a(0) \neq -F_a(1) = \calH(1)$ for $a \neq 1/2$. In particular, Theorem \ref{thm:Noether} does not hold in this case of a discontinuous solution.
\end{myremark}



\section{Hamiltonian equations}
\label{sec:Hamiltonian}
We further pursue the Noetherian formalism developed in Section \ref{sec:Noether} to interpret the conserved quantity associated with the translations as the Hamiltonian and construct the symplectic structure that formally associates the Hamiltonian with the vector field generating the translations in Section \ref{sec:sym}. The symplectic structure may well be degenerate and the Hamiltonian structure is understood only formally.  Section \ref{sec:cmfd} shows that the symplectic structure is automatically non-degenerate on finite-dimensional center manifolds as constructed in \cite{faye2016center}.

\subsection{Symplectic formalism}\label{sec:sym}
In the remainder of this paper, we restrict attention to a somewhat smaller class of Lagrangians, satisfying the following hypothesis.
\begin{myhyp}[L]
The Lagrangian $L$ is of the form
\[
L(u,v,n) = E(u,v) + \frac 1 2 S(u) \cdot n,
\]
where $E \in C^2(\R^d\times\R^d,\R)$.
\end{myhyp}
The corresponding Euler-Lagrange equation \eqref{eq:EL} is of the form
\begin{equation}
  \label{eq:EL_mechanical}
  - \partial_x \nabla_v E(u,u') + \nabla_u E(u,u') + DS(u)^T K*S(u) = 0.
\end{equation}

Define the function $\calH$ by
\begin{multline*}
  \calH(u) := - \rst{ \bigg( E(u,u') - \nabla_v E(u,u') \cdot u' + \frac 1 2 S(u) \cdot K * S(u) \bigg) }{ x = 0 } \\
  + \frac 1 2 \iint_Q S(u(x)) \cdot K(x-y) DS(u(y)) u'(y) - S(u(y)) \cdot K(x-y) DS(u(x)) u'(x) \d x \d y,
\end{multline*}
where  $Q := (-\infty,0) \times (0,\infty)$ is the upper left quadrant as in Theorem \ref{thm:Noether}.
We note that the symmetry group $G$ of \eqref{eq:EL_mechanical} contains the pure translations $\{\I_d\} \times \R$.
Let $\xi = (0 , 1) \in \frakg$, so that $u_\xi = u'$ and $\1_\xi = -1$.
Then $\calH(\phi_\tau u) = \calC_\xi(u)(\tau)$, where $\calC_\xi$ is the conserved quantity obtained in Theorem \ref{thm:Noether}.
An immediate consequence of Theorem\ref{thm:Noether} is then the following.
\begin{mycor}[Integral of motion]
  The quantity $\calH$ is conserved under the shift action $\phi_\tau$ on the solutions of \eqref{eq:EL_mechanical},
that is,
\[
\ddxe{\tau} \calH(\phi_\tau u) = 0, \qquad \text{for any $u$ solving \eqref{eq:EL_mechanical}}.
\]
\end{mycor}

Since $\calH$ corresponds to the translational symmetry of \eqref{eq:EL_mechanical} it is natural to think of $\calH$ as the
Hamiltonian of the system \eqref{eq:EL_mechanical}.
This raises the question whether the shift dynamics is Hamiltonian in a suitable sense.
It turns out that this is the case, at least formally, via a symplectic formalism.

\paragraph{Geometric setup.}
Before we introduce the relevant structures, let us briefly motivate the formalism we have in mind.
Say $\calM$ is a translation invariant set of functions $u : \R \to \R^d$ solving \eqref{eq:EL_mechanical},
and suppose that in some suitable topology $\calM$ comes equipped with the structure of a smooth manifold.
A presymplectic form $\omega$ is a closed $2$-form on $\calM$, i.e., for each $u \in \calM$ it defines a skew-symmetric bilinear form on $T_u \calM$.
If in addition $\omega$ is nondegenerate, it is called a symplectic form.
One can then define the Hamiltonian flow $\calX_\calH^\tau$ on the domain of definition of $D \calH$, by solving
\[
\partial_\tau \calX_\calH^\tau = V_\calH \circ \calX_\calH^\tau, \qquad \text{where} \quad - D \calH(u) = \omega_u\big( V_\calH(u) , \cdot \big).
\]
The aim is to find $\omega$ such that $V_\calH(u) = u'$, so that the shift symmetry $\calX_\calH^\tau(u) = \phi_\tau u$ is retrieved from the Hamiltonian flow of $\calH$.

The construction of a symplectic manifold $(\calM,\omega)$ consisting of solutions of \eqref{eq:EL_mechanical} is, in general, an open question.
We will address this question further in Section \ref{sec:cmfd}, where such a manifold is constructed for small amplitude solutions near (linear) center fixed points.
Instead of pursuing this matter further, in the current section we take $\calM$ to be a function space, say, $C^1_b(\R,\R^d)$.
Since solutions of \eqref{eq:EL_mechanical} only form a small subset of $C^1_b(\R,\R^d)$, 
in general it seems unlikely that here exists a symplectic structure defined on all of $C^1_b(\R,\R^d)$ such that the Hamiltonian flow of $\calH$ is the shift action $\phi_\tau$.
Instead we look for a presymplectic structure $\omega$ on $C^1_b(\R,\R^d)$ such that
\[
- D \calH(u) = \omega_u ( u' , \cdot ), \qquad \text{for any $u$ solving \eqref{eq:EL_mechanical}}.
\]
We will describe such a structure in the next paragraph.
If one obtains a manifold $\calM$ consisting of solutions of \eqref{eq:EL_mechanical}, the presymplectic structure is expected to restrict to a symplectic structure on $\calM$.
This is exemplified by the center manifold theory in Section \ref{sec:cmfd}.

\paragraph{Presymplectic structure.}
We will define an exact presymplectic structure on $C^1_b(\R,\R^d)$, hence we start with defining a $1$-form $\lambda$ on $C^1_b(\R,\R^d)$.
For $u \in C^1_b(\R,\R^d)$, define $\lambda^n_u : T_u C^1_b(\R,\R^d) \to \R$ by
\begin{equation}
  \label{eq:lambda_nonloc}
  \lambda^n_u(v) := \iint_Q v(x) \cdot DS(u(x))^T K(x-y) S(u(y)) \d x \d y, \qquad v \in T_u C^1_b(\R,\R^d) = C^1_b(\R,\R^d),
\end{equation}
with $Q =  (-\infty,0) \times (0,\infty)$.
In light of Hypothesis \ref{hyp:K}, the integral is convergent, cf.\ \eqref{eq:K_moment_Q_integrable}.
Note that $\lambda^n$ is $C^1$ in both $u$ and $v$.
Set $\omega^n := \d \lambda^n$, that is,
\[
\omega^n_u(v,w) := D_u [ \lambda^n_u(v) ](w) - D_u[ \lambda^n_u(w) ](v), \qquad v , w \in T_u C^1_b(\R,\R^d) = C^1_b(\R,\R^d).
\]
For readers familiar with exterior calculus we note that this expression follows from the usual invariant formula for exterior derivatives \cite{Lee2003},
which involves Lie derivatives along vector fields,
by considering vector fields obtained by constant extensions of the tangent vectors $v$ and $w$,
which is possible since our base manifold is a vector space.

A brief computation shows that
\[
\omega^n_u(v,w) = \iint_Q v(x) \cdot A_{xy}(u) w(y) - v(y) \cdot A_{yx}(u) w(x) \d x \d y,
\]
where
\[
A_{xy}(u) := DS(u(x))^T K(x-y) DS(u(y)).
\]
Thus we have defined an exact 2-form $\omega^n$ on $C^1_b(\R,\R^d)$.

We define another 1-form $\lambda^{\text{loc}}_u : T_u C^1_b(\R,\R^d) \to \R$ by
\begin{equation}
  \label{eq:lambda_loc}
  \lambda^{\text{loc}}_u(v) = - \rst{ \bigg( \nabla_v E(u,u') \cdot v \bigg) }{x=0}, \qquad v \in T_u C^1_b(\R,\R^d) = C^1_b(\R,\R^d).
\end{equation}
Then define a 2-form $\omega^{\text{loc}}_u := \d \lambda^{\text{loc}}_v$, that is,
\[
\omega^{\text{loc}}_u(v,w) := \rst{ \bigg( w\cdot \nabla^2_v E(u,u') v' - v \cdot \nabla^2_v E(u,u') w'
+ w \cdot \nabla_u \nabla_v E(u,u') v  - v \cdot \nabla_u \nabla_v E(u,u') w \bigg) }{ x=0 }.
\]

\begin{mydef}[Exact presymplectic structure]
Associated with \eqref{eq:EL_mechanical}, define the exact presymplectic structure $\omega$ on $C^1_b(\R,\R^d)$ as
\[
\omega := \d \lambda = \omega^n + \omega^{\text{loc}}, \qquad \text{where} \quad \lambda := \lambda^n + \lambda^{\text{loc}},
\]
with $\lambda^n$, $\lambda^{\text{loc}}$ given in \eqref{eq:lambda_nonloc}, \eqref{eq:lambda_loc}, respectively.
\end{mydef}

\begin{mylemma}[Formal Hamiltonian dynamics]
\label{lemma:Hamiltonian_time_shifts}
  Let $u \in C_b^R(\R,\R^d)$ be a solution of \eqref{eq:EL_mechanical}.
  For any $v \in T_u C^1_b(\R,\R^d) = C_b^1(\R,\R^d)$ we have
  \begin{equation}
    \label{eq:Hamiltonian_time_shifts}
    - D \calH(u) v = \omega_u( u' , v ).
  \end{equation}
\end{mylemma}
\begin{myproof}
Let $u \in C^R_b(\R,\R^d)$, $v \in  T_u C^1_b(\R,\R^d) = C^1_b(\R,\R^d)$ be arbitrary.
Note that
\begin{align*}
  \calH(u) &= - \rst{ \bigg( E(u,u') - \nabla_v E(u,u') \cdot u' + \frac 1 2 S(u) \cdot K * S(u) \bigg) }{ x = 0 } + \frac 1 2 \iint_Q S(u(x)) \cdot K(x-y) DS(u(y)) u'(y) \d x \d y \\
  &\quad - \frac 1 2 \iint_Q S(u(y)) \cdot K(x-y) DS(u(x)) u'(x) \d x \d y \\
  &= - \rst{ \bigg( E(u,u') + \frac 1 2 S(u) \cdot K * S(u) \bigg) }{ x = 0 } - \lambda^{\text{loc}}_u(u') \\
  &\quad - \frac 1 2 \lambda^n_u(u') + \frac 1 2 \iint_Q S(u(x)) \cdot K(x-y) DS(u(y)) u'(y) \d x \d y.
\end{align*}
To evaluate the last term, let $K_\epsilon \in C^\infty(\R,\calM(\R^d))$, $\epsilon > 0$ be such that $K_\epsilon$ satisfies hypothesis \ref{hyp:K}, and in addition,
\[
  \int_\R \left( 1 + |r| \right) | K_\epsilon'(r) | \d r < \infty,\qquad
  \lim_{\epsilon \to 0} \| K_\epsilon - K \|_{L^1(\R,\calM(\R^d))} = 0.
\]
Using integration by parts, the last term is transformed
\begin{align*}
          \frac 1 2 \iint_Q &S(u(x)) \cdot K_\epsilon(x-y) DS(u(y)) u'(y) \d x \d y \\
  &= \frac 1 2 \int_{x<0} \int_{y>0}  S(u(x)) \cdot K_\epsilon(x-y) \partial_y S(u(y)) \d x \d y \\
  &= \rst{- \frac 1 2 S(u) \cdot K_\epsilon*(\Theta_- S(u))}{x=0}  - \frac 1 2 \int_{x<0} \int_{y>0} S(u(x)) \cdot \big( \partial_y K_\epsilon(x-y) \big) S(u(y)) \d y \d x \\
  &= \rst{- \frac 1 2 S(u) \cdot K_\epsilon*(\Theta_- S(u))}{x=0} + \frac 1 2 \int_{y>0} \int_{x<0} S(u(x)) \cdot \big( \partial_x K_\epsilon(x-y) \big) S(u(y)) \d x \d y \\
  &= - \rst{\frac 1 2 S(u) \cdot K_\epsilon*(\Theta_- S(u))}{x=0} + \rst{ \frac 1 2 S(u) \cdot K_\epsilon*(\Theta_+ S(u))}{x=0}  \\
  &\quad - \frac 1 2 \int_{y>0} \int_{x<0} D S(u(x)) u'(x) \cdot K_\epsilon(x-y) S(u(y)) \d y \d x
  \end{align*}
where $\Theta_+(x) = 0$ for $x < 0$, $\Theta_+(x) = 1$ for $x \geq 0$, and $\Theta_- = 1 - \Theta_+$.
Letting $\epsilon \to 0$ in the final identity, we obtain
\begin{equation}
  \label{eq:Hamiltonian_boundary_term_transform}
          \frac 1 2 \iint_Q S(u(x)) \cdot K(x-y) DS(u(y)) u'(y) \d x \d y
          = \rst{\frac 1 2 \bigg( S(u) \cdot K*(\Theta_+ S(u)) - S(u) \cdot K*(\Theta_- S(u)) \bigg) }{x=0} - \frac 1 2 \lambda^n_u(u'),
\end{equation}
and thereby we obtain the fundamental relation between $\calH$ and $\lambda = \lambda^{\text{loc}} + \lambda^n$,
\[
\calH(u) = - \rst{ \bigg( E(u,u') + S(u) \cdot K*( \Theta_- S(u) ) \bigg) }{ x=0 } - \lambda_u(u').
\]
For the derivative we now obtain the expression
\begin{equation}
  \label{eq:derivative_HamSymp}
\begin{split}
  -D\calH(u)v &= \rst{  \bigg( \nabla_u E(u,u') \cdot v + \nabla_v E(u,u') \cdot v' \bigg)}{x=0} + D_u\big( \lambda^{\text{loc}}_u(u') \big) v \\
  &\quad + \rst{ \bigg( DS(u)^TK*( \Theta_- S(u) ) \cdot v + S(u) \cdot K*( \Theta_- DS(u) v ) \bigg) }{ x=0 } + D_u\big( \lambda^n_u(u') \big) v.
\end{split}
\end{equation}
Recalling the definition of $\omega^{\text{loc}} = \d \lambda^{\text{loc}}$, we have
\[
  D_u\big( \lambda^{\text{loc}}_u(u') \big) v = \omega^{\text{loc}}_u(u',v) + D_u\big( \lambda^{\text{loc}}_u(v) \big) u' + \lambda^{\text{loc}}_u(v').
\]
Now note that
\[
D_u\big( \lambda^{\text{loc}}_u(v) \big) u' + \lambda^{\text{loc}}_u(v') = - \rst{ \bigg( \partial_x \nabla_v E(u,u') \cdot v + \nabla_v E(u,u') \cdot v' \bigg) }{x=0},
\]
hence
\begin{multline}
  \label{eq:local_part_HamSymp}
     \rst{  \bigg( \nabla_u E(u,u') \cdot v + \nabla_v E(u,u') \cdot v' \bigg)}{x=0} + D_u\big( \lambda^{\text{loc}}_u(u') \big) v \\
   = \rst{  \bigg( - \partial_x \nabla_v E(u,u')  + \nabla_u E(u,u') \bigg)}{x=0} \cdot v(0) + \omega^{\text{loc}}_u(u',v).
\end{multline}

Computations similar to \eqref{eq:Hamiltonian_boundary_term_transform} show
\[
  D_u\big( \lambda^n_u(v) \big) u' + \lambda^n_u(v') = \rst{ \bigg( DS(u)^T K*( \Theta_+ S(u) ) \cdot v  -  S(u)\cdot K*(\Theta_- DS(u) v) \bigg) }{x=0}.
\]
Since $D_u\big( \lambda^n_u(u') \big) v = \omega^n_u(u',v) + D_u\big( \lambda^n_u(v) \big) u' + \lambda^n_u(v')$, it follows that
\begin{multline}
  \label{eq:nonlocal_part_HamSymp}
    \rst{ \bigg( DS(u)^TK*( \Theta_- S(u) ) \cdot v + S(u) \cdot K*( \Theta_- DS(u) v ) \bigg) }{ x=0 } + D_u\big( \lambda^n_u(u') \big) v \\
    \quad = \rst{ \bigg( DS(u)^T K*S(u) \bigg) }{x=0} \cdot v(0) + \omega^n_u(u',v).
\end{multline}
Combining \eqref{eq:derivative_HamSymp}, \eqref{eq:local_part_HamSymp}, and \eqref{eq:nonlocal_part_HamSymp}, we obtain
\[
  -D \calH(u) v = \omega_u(u',v)
+ \rst{ \bigg( - \partial_x \nabla_{v} E(u,u') + \nabla_u E(u,u') + DS(u)^T K*S(u) \bigg) }{x=0} \cdot v(0).
\]
From this the claimed result readily follows.
\end{myproof}

\begin{myremark}
In general, $\omega$ is degenerate as a $2$-form on $C^1_b(\R,\R^d)$.
To illustrate this, consider the following example.
Suppose the kernel has compact support, say, $\supp K \subset [-1,1]$.
Let $v$ have support $\supp v \subset \R \setmin [-1,1]$.
Then $\omega_u( v , \cdot ) = 0$, and this is true for any $u$. This degeneracy can be interpreted in several ways. First, when (formally) deriving \eqref{eq:Hamiltonian_time_shifts}, we notice that the relation holds with $u'\mapsto u'+\theta$ as long as $\supp\theta\subset \R \setmin [-1,1]$, such that the evolution equation is posed on $[-1,1]$, only, in the sense that the Hamiltonian formalism gives conditions on this interval, only. This is most apparent in the oversimplified case of $K(x)=\delta(x-1)+\delta(x+1)$, for which the equation decouples into a family of two-term recursions.  Clearly, bounded solutions can be found by solving an initial value problem with initial conditions given on $x\in [-1,1)$. The Hamiltonian equation on this interval can then be used to continue solutions beyond the interval. For equations including derivatives, for instance equations for breathers in nonlinear chains of oscillators, the resulting equations can still be viewed as evolution equations on an interval \cite{iooss}, which is however ill-posed. In our point of view, we understand dynamics as shifts on the set of bounded solutions without attempting to pose an initial-value problem. The Hamiltonian equation, formally derived here, would almost always be ill-posed. The absence of strong degeneracies for infinite-range kernels hints in this sense at the fact that no simple reduced formulation as a possibly ill-posed initial-value problem is avilable. 
\end{myremark}

\subsection[Local nondegeneracy of omega]{Local nondegeneracy of $\omega$ --- center manifold reduction}
\label{sec:cmfd}
We start this section with a brief summary of how the result from \cite{faye2016center} can be applied to equations of the form \eqref{eq:EL_mechanical},
to obtain a finite dimensional center manifold $\calM_0$, with the key property that all small solutions of \eqref{eq:EL_mechanical} are contained in $\calM_0$.
Since $\omega$ is typically degenerate on the ambient function space, our setup is somewhat different than the one considered in \cite{mielke2006hamiltonian}.
Nevertheless, we will show that the restriction of $\omega$ to $\calM_0$ is locally nondegenerate,
by explicit computations on the tangent space $T_0 \calM_0$.
Consequently, the center manifold $\calM_0$ naturally comes equipped with a symplectic structure $\omega$,
under which the induced Hamiltonian flow of $\calH$ is the shift flow $\phi_\tau$.

\paragraph{Construction of the center manifold.}
We will now summarize how the result from \cite{faye2016center} can be applied to equations of the form \eqref{eq:EL_mechanical}.
To do so, the equation \eqref{eq:EL_mechanical} needs to be recast in the form $\calT u + \calF(u) = 0$, posed on weighted Sobolev spaces allowing for exponential growth,
where $\calT u = u + J*u$ for some (not necessarily symmetric) $J$.
Given $\eta \in \R$, $1 \leq p \leq \infty$, $k \in \N \cup \{0\}$, let
\[
W^{k,p}_\eta(\R,\R^d) := \set{ u \in L^p_{\text{loc}}(\R,\R^d) }{ e^{\eta \langle x \rangle } \partial_x^j u(x) \in L^p(\R,\R^d), \; 0 \leq j \leq k }, \qquad L^p_\eta(\R,\R^d) := W^{0,p}_\eta(\R,\R^d),
\]
were $\langle x \rangle := \sqrt{ 1 + x^2 }$. From Morrey's inequality we have, for any $k \in \N \cup \{0\}$, $1 < p \leq \infty$, and $\eta \leq 0$, and suitable constants $C_{k,p,\eta}$ that
\begin{equation}
  \label{eq:Morrey_exponential}
  \| u \|_{C^k_b(\R,\R^d)} \leq C_{k,p,\eta} \sup_{\tau \in \R} \| \phi_\tau u \|_{W^{k+1,p}_\eta(\R,\R^d)}.
\end{equation}
We impose the following restriction upon \eqref{eq:EL_mechanical}.
\begin{myhyp}[C1]
\namedlabel{hyp:C1}{C1}
The zero function is a solution of \eqref{eq:EL_mechanical}, that is,
    \[
    \nabla_u E(0,0) + DS(0)^T \kappa_0 S(0) = 0, \qquad \text{with} \quad \kappa_0 := \int_\R K(r) \d r.
    \]
\end{myhyp}
In addition, we impose a local ellipticity condition. Define the principal symbol $\Sigma_p : \C \to \calM(\R^d)$ of \eqref{eq:EL_mechanical} at $u=0$ through
\begin{equation}
  \label{eq:principal_symbol}
  \Sigma_p(\nu) := - \nabla_v^2 E(0,0) \nu^2 + \nabla_u^2 E(0,0) + D^2S(0)^T \kappa_0 S(0).
\end{equation}
Let us also introduce the full symbol $\Sigma : \C \to \calM(\R^d)$ of \eqref{eq:EL_mechanical} at $u=0$, defined through
\begin{equation}
  \label{eq:full_symbol}
  \Sigma(\nu) := \Sigma_p(\nu) + DS(0)^T \widehat K(\nu) DS(0),\qquad \text{where} \quad \widehat K(\nu) := \int_\R e^{-\nu x} K(x) \d x.
\end{equation}
The exponential decay of $K$ ensures that $\widehat K(\nu) $ is analytic in a strip around $i \R$.
Note that the linearization of \eqref{eq:EL_mechanical} at $u=0$ takes the form
\[
\Sigma_p(\partial_x) u + \big( DS(0)^T K DS(0) \big) * u = 0,
\]
and $\Sigma(\nu) \widehat u(\nu) = 0$ is the (complex) Fourier transform of this equation.
\begin{myhyp}[C2]
\namedlabel{hyp:C2}{C2}
  We assume that $\det \Sigma_p(\nu) \neq 0$ for all $\nu \in i \R$.
  Furthermore, if $E(u,v)$ is not constant in $v$, we require that the map
  \[
  \lambda \mapsto |\lambda|^2 \| \Sigma_p(i \lambda)^{-1} \|, \qquad \text{with} \quad \lambda \in \R,
  \]
  is uniformly bounded.
\end{myhyp}
Our last hypothesis is concerned with smoothness of the pointwise evaluations and localization of the kernel.
\begin{myhyp}[C3]
\namedlabel{hyp:C3}{C3}
  \begin{enumerate}
  \item The kernel $K$ is exponentially localized, $K \in W^{1,1}_{\eta_0}(\R,\calM(\R^d))$ for some $\eta_0 > 0$.
  \item The function $S$ is of class $C^{k+R+2}$, $k \geq 2$.
  \item The function $E$ is of class $C^{k+2R+1}$ smooth, $k \geq 2$. 
  \end{enumerate}
  Here $R\in \{1,2\}$ is as in Definition \ref{def:solution_regularity}.
\end{myhyp}

Now now recast equation \eqref{eq:EL_mechanical} in the form $\calT u + \calF(u) = 0$, where
\[
\calT u := u + \Sigma_p(\partial_x)^{-1} \big( DS(0)^T K DS(0) \big) * u,
\]
and
\begin{equation}
  \label{eq:cmfd_F}
  \begin{split}
  \calF(u) := &\Sigma_p(\partial_x)^{-1} \bigg(  - \big( \nabla^2_v E(u,u') -  \nabla^2_v E(0,0) \big) u'' - \big( \nabla_u \nabla_v E(u,u') - \nabla_u \nabla_v E(0,0)  \big) u' \\
  &+ \big( \nabla_u E(u,u') - \nabla_u^2 E(0,0) u \big) + \big( DS(u)^T K* S(u) - D^2S(0)^T\kappa_0 S(0) u - (DS(0)^T K DS(0) )*u \big) \bigg).
  \end{split}
\end{equation}
Note that ellipticity of $\Sigma$ ensures that $\calF$ is locally $C^{k+R+1-r}$ on $W^{r,\infty}(\R,\R^d)$, for $R \leq r \leq k+R-1$. Our hypotheses ensure that the maps $\calT$ and $\calF$ satisfy Hypotheses (H1) and (H2) from \cite{faye2016center}. As usual, the global center manifold is only defined for sufficiently small nonlinearities $\calF$.
A local center manifold is then defined by modifying the original nonlinearity $\calF$ outside of a small neighborhood of $0$. We want to point out that in this setting, this modification can be made without destroying the variational structure. Indeed, define smooth cutoff functions $\chi_\epsilon : \R \to \R$ by setting
\[
\chi_1(x) :=
\begin{cases}
  1 \quad &\text{if } |x| \leq 1, \\
  0 \quad &\text{if } |x| \geq 2,
\end{cases}
\]
and
\[
\chi_\epsilon(x) := \chi_1(x/\epsilon).
\]
Then replace $E$ and $S$ in \eqref{eq:EL_mechanical} by $E^\epsilon(u,v) := \chi_\epsilon(|(u,v)|) E(u,v)$ and $S^\epsilon(u) := \chi_\epsilon(u) S(u)$, respectively.
Now let $\calF^\epsilon$ be as in \eqref{eq:cmfd_F}, with $E$ and $S$ replaced by $E^\epsilon$ and $S^\epsilon$.
For sufficiently small $\epsilon > 0$, such a modified nonlinearity will satisfy the smallness assumptions needed in the center manifold theorem.

The construction of the center manifold itself can then be summarized as follows.
First, choose $r$ such that $R \leq r \leq k+R-1$.
We will expand the center manifold as a subspace of $W^{r,2}_{-\delta}(\R,\R^d)$.
One sets $\calE_0 := \ker(\calT) \subset W^{r,2}_{-\eta}(\R,\R^d)$ for $0 < \eta < \eta_0$.
To avoid trivial situations, we will assume $\calE_0 \neq \{0\}$.
It follows from the exponentially localized behavior of $K$ that $\calE_0$ is finite dimensional, independent of $\eta$ for sufficiently small $\eta$,
and the dimension is given by the algebraic count
of the number of solutions (counting multiplicity) of the characteristic equation
\begin{equation}
  \label{eq:characteristic}
  \det \Sigma(\nu) = 0, \qquad \text{where} \quad \nu = i \ell \in i \R.
\end{equation}
Here $\Sigma$ is the full symbol defined in \eqref{eq:full_symbol}.

One has the freedom to choose a projection $\calQ : W^{r,2}_{-\eta}(\R,\R^d) \to W^{r,2}_{-\eta}(\R,\R^d)$ onto $\calE_0$. 
This projection should be equivariant with respect to the shift map $\phi_\tau$.
Furthermore, the projection $\calQ$ has to be constructed in a way which is independent of the weight $\eta > 0$.
We remark that such a projection always exists, and the particular choice does not influence the abstract results in this section and the next,
however from a computational viewpoint certain choices can be favorable.
Then there exists a $\delta >0$ and a $C^{k+R+1-r}$ smooth map
\[
\Psi : \calE_0 \to \ker \calQ \subset W^{r,2}_{-\delta}(\R,\R^d),\qquad \Psi(0) = 0,\  D \Psi(0) = 0,
\]
such that the graph
\[
\calM_0 := \set{ u_0 + \Psi(u_0) }{ u_0 \in \calE_0 }
\]
consists precisely of the solutions $u \in W^{r,2}_{-\delta}(\R,\R^d)$ to the modified equation ${\calT u + \calF^\epsilon(u) = 0}$.
Furthermore, $\Psi$ commutes with the shift map, $\Psi \circ \phi_\tau = \phi_\tau \Psi$.
Any solution $u \in W^{r,2}_{-\delta}(\R,\R^d)$ of equation \eqref{eq:EL_mechanical},
which is sufficiently small, in the sense that $\| u \|_{C^{R-1}_b(\R,\R^d)} < \epsilon$, is contained in $\calM_0$,
simply because such a $u$ will solve the modified equation ${\calT u + \calF^\epsilon(u) = 0}$.
Conversely, in light of inequality \eqref{eq:Morrey_exponential}, if $u \in \calM_0$ whose entire orbit 
\[
\gamma(u) := \bigcup_{\tau \in \R} \left\{ \phi_\tau u \right\}
\]
is sufficiently small in $W^{r,2}_{-\eta}(\R,\R^d)$, it follows that $u$ solves \eqref{eq:EL_mechanical}.

\paragraph{Symplectic structure on the center manifold.}
  
We next construct an explicit basis for $\calE_0$ in which we establish nondegeneracy of $\omega_0$ restricted to $\calE_0$. Let $\nu_1,\dots,\nu_n$ denote the roots of the characteristic equation \eqref{eq:characteristic}, which lie in the upper half plane, that is $\nu_j = i \ell_j$ with $\ell_j \geq 0$.
Given $1 \leq j \leq n$, let $\gamma_j := \dim \ker \Sigma(\nu_j)$, and choose a real basis $\e_{j,1},\dots,\e_{j,\gamma_j}$ for $\ker \Sigma(\nu_j)$.
We now impose the following restriction upon $\Sigma$.
\begin{myhyp}[C4]
  \namedlabel{hyp:C4}{C4}
For each $1 \leq j \leq n$ and $1 \leq k \leq \gamma_j$,
there exists a $n_{j,k} \in \N$ such that
\[
\Sigma( \nu ) \e_{j,k} = \frac{1}{n_{j,k}!} \Sigma^{(n_{j,k})}(\nu_j) ( \nu - \nu_j )^{n_{j,k}} \e_{j,k} + \calO\left( |\nu-\nu_j|^{n_{j,k}+1}\right) \qquad \text{as} \quad \nu \to \nu_j.
\]
\end{myhyp}

Given $1 \leq k \leq \gamma_j$ and $0 \leq l \leq n_{j,k}-1$, define
\[
\psi_{j,k,l}(x) := x^l \left( \cos\big( \ell_j x \big) + \sin\big( \ell_j x \big) \right) \e_{j,k}.
\]
Then define the conjugate vector $\psi^*_{j,k,l}$ by
\[
\psi^*_{j,k,l}(x) :=
\begin{cases}
  x^{n_{j,k} - l - 1} \left( \cos\big( \ell_j x \big) + \sin\big( \ell_j x \big) \right) \e_{j,k} \quad & \text{if $n_{j,k}$ is even}, \\
  x^{n_{j,k} - l - 1} \left( \cos\big( \ell_j x \big) - \sin\big( \ell_j x \big) \right) \e_{j,k} \quad & \text{if $n_{j,k}$ is odd}.
\end{cases}
\]

\begin{mylemma}
  The system of vectors
\[
B := \set{  \psi_{j,k,l} , \;  \psi^*_{j,k,l} }{ 1 \leq j \leq n, \; 1 \leq k \leq \gamma_j, \; 0 \leq l \leq n_{j,k} - 1 }
\]
is a basis for $\calE_0$.
\end{mylemma}
\begin{myproof}
Note that by symmetry of $K$, if $\nu_j = 0$ for some $1 \leq j \leq n$, then $n_{j,k}$ is even, for each $1 \leq k \leq \gamma_j$.
Hence the vectors always come in pairs $( \psi_{j,k,l} , \psi^*_{j,k,l} )$, with $\psi_{j,k,l}$ and $\psi^*_{j,k,l}$ linearly independent.
Taking the distributional complex Fourier transform of $\calT \psi_{j,k,l}$, we get
\[
\Sigma_p( \nu_j ) \widehat{ \calT \psi_{j,k,l} } = 
\begin{cases}
  \Sigma^{(l)}(\nu_j) \e_{j,k} \quad & \text{if $l$ is even}, \\
  - i \Sigma^{(l)}(\nu_j) \e_{j,k} \quad & \text{if $l$ is odd}.
\end{cases}
\]
Hence, by Hypothesis \ref{hyp:C4}, we have $\psi_{j,k,l} \in \calE_0$, and, similarly, we find $\psi^*_{j,k,l} \in \calE_0$.
We note that $\alpha_j := n_{j,1} + \cdots + n_{j,k}$ is the algebraic multiplicity of the root $\nu_j$ of the characteristic equation \eqref{eq:characteristic}.
Since $\dim \calE_0$ is given by the number of roots of \eqref{eq:characteristic}, counting multiplicity,
we conclude that the system of vectors $B$ forms a basis for $\calE_0$.
\end{myproof}

The vectors $\psi_{j,k,l}$, $\psi^*_{j,k,l}$ turn out to be each other's symplectic dual.
In preparation to proving this fact, we have the following useful lemma.
\begin{mylemma}
  \label{lemma:omega_covariance_computation}
  Fix $p, q \in \N \cup \{0\}$.
  Suppose $f : \R^2 \to \R$ is of the form
  \[
  f(x,y) = g(x-y) + h(x+y), \qquad \text{where} \quad g , h \in L^2_{-\eta}(\R,\R),
  \]
with $g$ and $h$  of the form
  \[
  \begin{array}{l l l}
    g(r) = -g(-r), \quad & h(r) = h(-r), \qquad &\text{if $p+q$ is even}, \\
    g(r) = g(-r), \quad & h(r) = -h(-r), \qquad &\text{if $p+q$ is odd}.
  \end{array}
  \]
  Then
  \[
    \iint_Q \bigg( x^p y^q f(x,y) - x^q y^p f(y,x) \bigg) K(x-y) \d x \d y = C_{p,q} \int_\R r^{p+q+1} g(r) K(r) \d r,
  \]
  where
  \[
  C_{p,q} := \frac{ (-1)^q }{ 2^{p+q} } \sum_{\substack{ 0 \leq m \leq p \\ 0 \leq n \leq q}} { p \choose m } { q \choose n } \frac{ \sigma_{m,n} }{ m+n+1 },
  \]
  with
  \[
  \sigma_{m,n} :=
  \begin{cases}
    (-1)^{m+1} \quad &\text{if $m+n$ is even}, \\
    0 \quad &\text{if $m+n$ is odd}.
  \end{cases}
  \]
  In particular,
  \[
  \iint_Q \bigg( x^p y^q - x^q y^p \bigg) K(x-y) \d x \d y = C_{p,q} \widehat K^{(p+q+1)}(0).
  \]
\end{mylemma}
\begin{myproof}
  Let $s : \R^2\to \R^2$ be the reflection in the anti-diagonal, that is, $s(x,y) = -(y,x)$.
  Let $s^*$ denote the induced action on functions $F : \R^2 \to X$, defined through the pullback relation
  \[
  s^* F = F \circ s.
  \]
  Then
  \[
  s^{-1}(Q) = Q, \qquad s^* \big( x^q y^p f(x,y) \big) = - x^p y^q g(x-y) + x^p y^q h(x+y).
  \]
  Therefore,
  \begin{equation}
    \label{eq:2}
    \iint_Q \bigg( x^p y^q f(x,y) - x^q y^p f(y,x) \bigg) K(x-y) \d x \d y = 2 \iint_Q x^p y^q g(x-y) K(x-y) \d x \d y.
  \end{equation}
  We now reparameterize $Q$ by $    (x,y) = (( r + s ) / 2, ( s - r ) / 2)$ where  $r < 0, r < s < -r$.
  and compute,
  \begin{align*}
    \iint_Q &x^p y^q g(x-y) K(x-y) \d x \d y = \frac{1}{2^{p+q+1}} \int_{r < 0} \int_{|s| < |r|} (s + r)^p (s-r)^qg(r) K(r) \d s \d r \\
    &= \frac{1}{2^{p+q+1}} \int_{r < 0} \int_{|s| < |r|}  \bigg( \sum_{m=0}^p { p \choose m } r^{p-m} s^m \bigg) \bigg( \sum_{n=0}^p { q \choose n } (-r)^{q-n} s^n\bigg)  g(r) K(r) \d s \d r \\
    &= \frac{(-1)^q}{2^{p+q+1}} \sum_{\substack{ 0 \leq m \leq p \\ 0 \leq n \leq q}} { p \choose m } { q \choose n } (-1)^{n} 
    \int_{r < 0} \left( \int_{|s| < |r|} s^{m+n} \d s \right) r^{p+q-m-n} f(r) K(r) \d r \\
    &= \frac{(-1)^q}{2^{p+q}} \sum_{\substack{ 0 \leq m \leq p \\ 0 \leq n \leq q}} { p \choose m } { q \choose n } \frac{ \sigma_{m,n} }{ m+n+1 } \int_{r < 0} r^{p+q+1} g(r) K(r) \d r \\
    &= \frac 1 2 C_{p,q} \int_\R r^{p+q+1} g(r) K(r) \d r.
  \end{align*}
  Combined with \eqref{eq:2}, the claim follows.
\end{myproof}

For future reference, we now collect the hypothesis from this section.
\begin{myhyp}[C$\omega$]
  \namedlabel{hyp:Comega}{C$\omega$}
  Assume Hypotheses \ref{hyp:C1}, \ref{hyp:C2}, \ref{hyp:C3}, and \ref{hyp:C4} are satisfied.
\end{myhyp}

We are now prepared to show that the center manifold is locally symplectic.
\begin{mythm}[Symplectic center manifold]
\label{thm:symplectic_center_manifold}
  Assume Hypothesis \ref{hyp:Comega}and $R \leq r \leq k+R-1$.
  Then the following holds.
  \begin{enumerate}
  \item The center manifold $\calM_0 \subset W^{r,2}_{-\delta}(\R,\R^d)$ is $C^{k+R+1-r}$ smooth, and the flow $\phi_\tau$ is $C^{k+R-1-r}$ smooth on $\calM_0$.
  \item The Hamiltonian $\calH$ and presymplectic form $\omega = \d \lambda$ have $C^{k+R+1-r}$ smooth restrictions to $\calM_0$.
  \item The presymplectic form $\omega$ is locally nondegenerate, hence symplectic.
  \end{enumerate}
  Consequently, there exists $\delta > 0$, $\epsilon > 0$ so that if $u \in W^{r,2}_{-\delta}(\R,\R^d)$ solves \eqref{eq:EL_mechanical},
  and ${\|u(x)\|_{C^{R-1}_b(\R,\R^d)} < \epsilon}$, then $u \in \calM_0$.
  Furthermore, there exists an open neighborhood $0 \in O \subset \calM_0$ on which the Hamiltonian flow $\calX_\calH^\tau$, given by
  \[
  \partial_\tau \calX_\calH^\tau = V_\calH \circ \calX_\calH^\tau, \qquad \text{where} \quad - D \calH(u) = \omega_u\big( V_\calH(u) , \cdot \big),
  \]
  is defined for small $\tau$.
  Moreover, if $R \leq r \leq k+R-2$, we have $V_\calH(u) = u'$ and $\calX_\calH^\tau(u) = \phi_\tau u$ for all $u \in O$.
\end{mythm}
\begin{myproof}
We unfold a scale of center manifolds 
\[
\calM_0^r \subset W^{r,2}_{-\delta}(\R,\R^d), \qquad \text{for} \quad R \leq r \leq k +R - 1.
\]
Then $\calM_0^r$ is a $C^{k+R+1-r}$ smooth manifold.
In particular, the manifold $\calM_0^{k+R-1}$ is $C^2$ smooth, but the flow $\phi_\tau$ is only $C^0$.
Of course, given $u \in \calM_0^{k+R-1}$, the curve $\tau \mapsto \phi_\tau u$ is $C^{k+R-1-r}$ as a map from $\R$ into $W^{r,2}_{-\delta}(\R,\R^d)$.
Now we let
\[
\iota^{r,k+R-1} : W^{k+R-1}_{-\delta}(\R,\R^d) \hookrightarrow W^{r,2}_{-\delta}(\R,\R^d)
\]
 denote the natural inclusion between Sobolev spaces.
Then, it follows that $\phi_\tau$ is a $C^{k+R-1-r}$ smooth map in the inherited topology on 
\[
\iota^{r,k+R-1}( \calM_0^{k+1} ) \subset W^{r,2}_{-\delta}(\R,\R^d).
\]
It follows from the local uniqueness of the center manifold that 
\[
\iota^{r,k+R-1}( \calM_0^{k+R+r-1} ) = \calM_0^r,
\]
hence in particular it is $C^{k+R+1-r}$ smooth.
In particular, on $\calM_0 \subset W^{R,2}_{-\delta}(\R,\R^d)$ the flow $\phi_\tau$ will always be at least $C^1$ smooth.

Concerning the extensions of $\calH$ and $\omega = \d \lambda$.
we first recall that the modified equation $\calT u + \calF^\epsilon(u) = 0$ is again of the form \eqref{eq:EL_mechanical}.
Hence associated with the modified equation we obtain Hamiltonians $\calH^\epsilon$, and presymplectic structures $\omega^\epsilon = \d \lambda^\epsilon$.
We have $\calH(u) = \calH^\epsilon(u)$, and $\omega_u = \omega^\epsilon_u$, whenever $u$ satisfies $\|u(x)\|_{C^{R-1}_b(\R,\R^d)} < \epsilon$.
Since the center manifold only contains relevant information about \eqref{eq:EL_mechanical} for small amplitude solutions,
we may thus without loss of generality replace $\calH$ and $\omega = \d \lambda$ with their localized counterparts.

Note that in the definition of $\lambda$ and $\omega$ we have worked with spaces of uniformly bounded functions,
but the center manifold $\calM_0$ is naturally unfolded in the exponentially weighted space $W^{r,2}_{-\eta}(\R,\R^d)$.
Repeating the definition of $\lambda^n$ and $\omega^n$, but this time using functions from $W^{r,2}_{-\delta}(\R,\R^d)$,
one readily finds that the constructions are still valid, provided that the maps
\begin{align*}
  (x,y) &\mapsto S(u(x)) \cdot K(x-y) DS(u(y)) v(y), \\
  (x,y) &\mapsto v(x) \cdot DS(u(x))^T K(x-y) DS(u(y)) w(y)
\end{align*}
are integrable over $(x,y) \in Q$, for any $u , v , w \in W^{r,2}_{-\delta}(\R,\R^d)$.
This is clearly the case for the localized nonlinearities.
Consequently, $\lambda$ and $\omega$ then extend to $C^{k+R+1-r}$ structures on $W^{r,2}_{-\delta}(\R,\R^d)$.
Similarly, we find that $\calH$ extends to a $C^{k+R+1-r}$ smooth function on $W^{r,2}_{-\delta}(\R,\R^d)$.

  Since $\omega_u$ depends $C^{k+R+1-r}$ smoothly on the base point $u \in \calM_0$,
  it suffices to check whether $\omega_0$ is nondegenerate on $T_0 \calM_0 = \calE_0$.
  Suppose first that $n_{j,k}$ is even.
  Then, with the help of Lemma \ref{lemma:omega_covariance_computation}, we have
  \begin{align*}
    \omega^n_0\big( \psi_{j,k,l} , \psi^*_{j,k,l} \big) &= C_{l,n_{j,k}-l-1} \e_{j,k} \cdot DS(0)^T \left( \int_\R r^{n_{j,k}} \cos( \ell_j r ) K(r) \d r \right) DS(0) \e_{j,k} \\
    &= \frac 1 2 C_{l,n_{j,k}-l-1} \e_{j,k} \cdot DS(0)^T \left( \widehat K^{(n_{j,k})}(\nu_j) + \widehat K^{(n_{j,k})}(-\nu_j) \right) DS(0) \e_{j,k} \\
    &= C_{l,n_{j,k}-l-1} \e_{j,k} \cdot DS(0)^T \widehat K^{(n_{j,k})}(\nu_j) DS(0) \e_{j,k} \\
    &= C_{l,n_{j,k}-l-1} \e_{j,k} \cdot \bigg( \Sigma^{(n_{j,k})}(\nu_j) - \Sigma_p^{(n_{j,k})}(\nu_j) \bigg) \e_{j,k},
  \end{align*}
  where we used that $\widehat K(\nu) = \widehat K(-\nu)$.
  One observes that
  \[
  \omega^{\text{loc}}_0\big( \psi_{j,k,l} , \psi^*_{j,k,l} \big) = \Sigma_p^{(n_{j,k})}(\nu_j)  = 0 \qquad \text{whenever} \quad n_{j,k} \geq 3.
  \]
  For $n_{j,k} = 2$, a case-by-case examination shows that
  \[
  \omega^{\text{loc}}_0\big(\psi_{j,k,l} , \psi^*_{j,k,l} \big) = C_{l,n_{j,k}-l-1} \e_{j,k} \cdot \Sigma_p^{(n_{j,k})}(\nu_j) \e_{j,k},
  \]
  Consequently, by Hypothesis \ref{hyp:C4},
  \[
  \omega_0\big(\psi_{j,k,l} , \psi^*_{j,k,l} \big) = C_{l,n_{j,k}-l-1} \e_{j,k} \cdot \Sigma^{(n_{j,k})}(\nu_j) \e_{j,k} \neq 0.
  \]
  Similarly, when $n_{j,k}$ is odd, we have
  \begin{align*}
    \omega^n_0\big( \psi_{j,k,l} , \psi^*_{j,k,l} \big) &= C_{l,n_{j,k}-l-1} \e_{j,k} \cdot DS(0)^T \left( \int_\R r^{n_{j,k}} \sin( \ell_j r ) K(r) \d r \right) DS(0) \e_{j,k} \\
    &= - i C_{l,n_{j,k}-l-1}  \e_{j,k} \cdot \bigg( \Sigma^{(n_{j,k})}(\nu_j) - \Sigma_p^{(n_{j,k})}(\nu_j) \bigg) \e_{j,k}.,
  \end{align*}
  and
  \[
  \omega^{\text{loc}}_0\big(\psi_{j,k,l} , \psi^*_{j,k,l} \big) = - i C_{l,n_{j,k}-l-1} \e_{j,k} \cdot \Sigma_p^{(n_{j,k})}(\nu_j) \e_{j,k},
  \]
  so that, again,
  \[
  \omega_0\big(\psi_{j,k,l} , \psi^*_{j,k,l} \big) = - i C_{l,n_{j,k}-l-1} \e_{j,k} \cdot \Sigma^{(n_{j,k})}(\nu_j) \e_{j,k} \neq 0.
  \]
  This shows that $\omega_0$ is nondegenerate.

We obtain an open neighborhood $0 \in O \subset W^{r,2}_{-\eta}(\R,\R^d)$ on which $\omega$ is nondegenerate.
We then obtain a Hamiltonian flow $\calX_\calH^\tau$ on $O$, induced by $\calH$ via the symplectic form $\omega$.
If $R \leq r \leq k+R-2$, the flow $\phi_\tau$ on $\calM_0$ is differentiable,
so that $u' \in T_u \calM_0$ for any $u \in O$.
Consequently, in light of Lemma \ref{lemma:Hamiltonian_time_shifts}, we have $V_\calH(u) = u'$ and $\calX_\calH^\tau(u) = \phi_\tau u$ for all $u \in O$.
\end{myproof}

As a cautionary note, we want to point out here that the vectors $\psi_{j,k,l}$, $\psi^*_{j,k,l}$ are, in general, not in Darboux normal form, 
as the following example illustrates.
\begin{myexample}
  Let $K : \R \to \R$ be a kernel for which
  \[
  \Sigma^{(r)}(0) = 0, \quad \Sigma^{(4)}(0) \neq 0, \qquad \text{for} \quad 0 \leq r \leq 3,
  \]
  and suppose
  \[
  \widehat K^{(6)}(0) \neq 0.
  \]
Basis vectors for $\calE_0$ are then
\[
\psi_{l}(x) := x^l, \quad \psi^*_l(x) := x^{3-l}, \qquad \text{for} \quad 0 \leq l \leq 3,
\]
One then readily finds
\[
\omega_0\big( \psi_2 , \psi_3 \big) = \omega_0\big( \psi^*_1,\psi^*_0 \big) = \omega_0\big( \psi_2,\psi^*_0 \big) = \omega_0\big( \psi^*_1 , \psi_3 \big) = C_{2,3} \widehat K^{(6)}(0) \neq 0.
\]
\end{myexample}

\paragraph{Symplectic structure on a parameter dependent center manifold.}
Suppose that $K = K_\mu$, $E= E_\mu$, and $S = S_\mu$ depend on a parameters $\mu = (\mu_1,\dots,\mu_s) \in (-\epsilon',\epsilon')^s$.
We then impose the following conditions.
\begin{myhyp}[C$\mu$]
\namedlabel{hyp:Cmu}{C$\mu$}
  \begin{enumerate}
  \item In $\mu = 0$, hypothesis \ref{hyp:C1} is satisfied, that is,
    \[
    \nabla_u E_0(0,0) + DS_0(0)^T \kappa_0 S_0(0) = 0.
    \]
  \item Hypothesis \ref{hyp:C2} holds, uniformly in $\mu$.
    More precisely, let $\Sigma_p(\nu ; \mu)$ denote the principal symbol corresponding to the parameter value $\mu$.
    We demand that $\det \Sigma_p(\nu) \neq 0$ for all $\nu \in i \R$ and $\mu \in (-\epsilon',\epsilon')^s$.
    Furthermore, if for some parameter $\mu \in (-\epsilon',\epsilon')^s$ the map $E_\mu(u,v)$ is not constant in $v$, we require that the map
  \[
  (\lambda,\mu) \mapsto |\lambda|^2 \| \Sigma_p(i \lambda ; \mu)^{-1} \|, \qquad \text{with} \quad \lambda \in \R, \quad \mu \in (-\epsilon',\epsilon')^s,
  \]
  is uniformly bounded.
  \item Hypothesis \ref{hyp:C3} holds, with all estimates uniform in $\mu$.
  \item Hypothesis \ref{hyp:C4} holds for $\mu = 0$.
  \item The maps $K_\mu$, $E_\mu$, and $S_\mu$ depend $C^{k+4}$ smoothly on the parameter $\mu$.
  \end{enumerate}
\end{myhyp}

As before, one then localizes the nonlinearities $E_\mu$ and $S_\mu$,
and recasts the modified equation into the form $\calT_\mu u + \calF_\mu^\epsilon(u) = 0$, where
\[
\calT_\mu u = u + \Sigma_p(\partial_x ; \nu)^{-1} \big( DS_\mu(0)^T K DS_\mu(0) \big) * u.
\]
Then let $\calE_0 := \ker \calT_0$ and choose a projection $\calQ : W^{R,2}_{-\eta}(\R,\R) \to W^{R,2}_{-\eta}(\R,\R)$ onto $\calE_0$, as before.
After shrinking $\epsilon'$ if necessary, one then obtains a map
\[
\Psi : \calE_0 \times (-\epsilon',\epsilon')^s \to \ker \calQ,
\]
with $\Psi(0,0) = D_u \Psi(0,0) = 0$, such that
\[
\calM_0(\mu) := \set{ u_0 + \Psi(u_0,\gamma,\mu) }{ u_0 \in \calE_0 }, \qquad \text{with} \quad \mu \in (-\epsilon',\epsilon')^s,
\]
consists precisely of the solutions $u \in W^{R,2}_{-\delta}(\R,\R^d)$ of the modified equation $\calT_\mu u + \calF_\mu^\epsilon(u) = 0$.
Conversely, if $u \in W^{R,2}_{-\delta}(\R,\R^d)$ is a solution of \eqref{eq:EL_mechanical} corresponding to a parameter value $|\mu| < \epsilon'$,
and $\|u\|_{C^{R-1}_b} < \epsilon$, it follows that $u \in \calM_0(\mu)$.
Thus, we have obtained a local reduction of the parameter dependent equation \eqref{eq:EL_mechanical}.
Furthermore, the results from Theorem \ref{thm:symplectic_center_manifold} extend to the parameter dependent setting.
Hence, we obtain Hamiltonians $\calH_\mu$ and local symplectic forms $\omega_\mu = \d \lambda_\mu$ on $\calM_0(\mu)$, which are $C^k$ in $u \in \calM_0(\mu)$
as well as $C^k$ smooth in the parameter $\mu$.
The Hamiltonian flow of any such $\calH_\mu$ coincides with shift action $\phi_\tau$ on $\calM_0(\mu)$.



\section{Dissipative equations}
\label{sec:grad_like}
We adapt our formalism to traveling-wave equations, where the variational structure induces a gradient-like structure. Section \ref{sec:gl} contains the general formalism and Section \ref{sec:glcmfd} contains adaptations to local center manifolds. 

\subsection{Gradient-like behavior} \label{sec:gl}
Let $\Gamma \in C^0_b(\R^d,\calM(\R^d))$ take its values in the cone of positive definite $d \times d$-matrices.
Formally consider the equation $u_t = - \Gamma(u)^{-1} \nabla_{L^2} \calA_{-\infty}^\infty(u)$,
which is a gradient flow on a Hilbert manifold $\calM$ modeled over $L^2(\R,\R^d)$, with Riemannian metric given by
\[
g_u(v,w) := \langle v , \Gamma(u) w \rangle_{L^2(\R,\R^d)}, \qquad \text{where} \quad v, w \in T_u \calM.
\]
We make the traveling wave ansatz $u(t,x) = u(\xi)$, $\xi = x - c t$, with $c \neq 0$.
This leads to an equation of the form
\begin{equation}
  \label{eq:grad_like}
  - \partial_x \nabla_v E(u,u') + \nabla_u E(u,u') + DS(u)^T K*S(u)) = c \Gamma(u) u',
\end{equation}
We will consider $C^R$ smooth solutions, where $R = 1$ if $\nabla_v L = 0$, or $R = 2$ otherwise.

Since \eqref{eq:grad_like} is obtained by considering a gradient flow in a comoving frame,
we expect \eqref{eq:grad_like} to be a gradient-like system.
Indeed, this turns out to be the case.
Define the Lyapunov function $\calL$ by
\begin{multline*}
  \calL(u) := - \rst{ \bigg( E(u,u') - \nabla_v E(u,u') \cdot u' + \frac 1 2 S(u) \cdot K * S(u) \bigg) }{ x = 0 } \\
  + \frac 1 2 \iint_Q S(u(x)) \cdot K(x-y) DS(u(y)) u'(y) - S(u(y)) \cdot K(x-y) DS(u(x)) u'(x) \d x \d y.
\end{multline*}
We then have the following result.
\begin{mylemma}[Lyapunov function]
\label{lemma:Lyapunov}
  The quantity $\calL$ is monotone under the shift action $\phi_\tau$ on solutions of \eqref{eq:grad_like}, 
  \[
  \calL(\phi_b u) - \calL(\phi_a u) = - c \int_a^b \Gamma(u) u' \cdot u' \d x, \qquad \text{for any $u$ solving \eqref{eq:grad_like},
  and $-\infty < a \leq b < \infty$}.
  \]
\end{mylemma}
\begin{myproof}
We note that the symmetry group $G$ of \eqref{eq:grad_like} contains the pure translations $\{\I_d\} \times \R$.
Let $\xi = (0 , 1) \in \frakg$, so that $u_\xi = u'$ and $\1_\xi = -1$.
Then $\calL(\phi_\tau u) = \calC_\xi(u)(\tau)$, where $\calC_\xi$ is the conserved quantity corresponding to $\xi$.
Inspecting the proof of Theorem \ref{thm:Noether}, we note that the Euler-Lagrange equation \eqref{eq:EL} is used  in \eqref{eq:noether_step}, only.
If $u$ in turn solves \eqref{eq:grad_like}, we obtain
\begin{multline*}
    - c \int_a^b \Gamma(u) u' \cdot u_\xi \d x = \eval{ L\big(x,u(x),u'(x),K*S(u)(x)\big) \1_\xi }{x=a}{b} + \eval{ \nabla_v L(x,u,u',K*S(u)) \cdot u_\xi }{x=a}{b} \\
    + \int_a^b \nabla_n L\big(x,u,u',K*S(u)\big) \cdot K*DS(u) u_\xi
    - DS(u)^T K* \nabla_n L\big(x,u,u',K*S(u)\big) \cdot u_\xi \d x.
\end{multline*}
The remainder of the proof of Theorem \ref{thm:Noether} transfers without changes,
resulting in the desired expression
\[
\calC_\xi(u)(b) - \calC_\xi(u)(a) = - \int_a^b \Gamma(u) u' \cdot u_\xi \d x.
\]
\end{myproof}

We now impose the following additional conditions upon the Lagrangian $L$.
\begin{myhyp}[G]
  \namedlabel{hyp:G}{G}
\begin{enumerate}
\item The functions $E$ and $S$ are $C^{R+1}$ smooth.
\item For all $u$, $v$, either $\nabla_v E(u,v) = 0$, or $\nabla^2_v E(u,v)$ is invertible.
Here $\nabla^2_v L$ denotes the Hessian with respect to the variable $v$.
\end{enumerate}
\end{myhyp}

\begin{mydef}[$\alpha$- and $\omega$-limit sets]
\label{def:limit_sets}
  Given $N \subset C_{\text{loc}}^R(\R,\R^d)$,
we define $\alpha(N)$ to be the $\alpha$-limit set of $N$ with respect to the shift dynamics on $C_{\text{loc}}^R(\R,\R^d)$.
Thus, $\alpha(N)$ consists of the accumulation points of $\set{ \phi_\tau u }{ \tau < 0,\; u \in N }$.
Similarly, the $\omega$-limit set $\omega(N)$ is defined as the set of accumulation points of  $\set{ \phi_\tau u }{ \tau > 0,\; u \in N }$.
\end{mydef}

\begin{mythm}[Gradient-like behavior]
\label{thm:grad_like}
  Assume $L$ satisfies Hypothesis \ref{hyp:G}.
  Let $u \in C^R_b(\R,\R^d)$ be a solution of \eqref{eq:grad_like}, bounded in the $C^{R-1}_b(\R,\R^d)$-norm.
  Then both $\alpha(u)$ and $\omega(u)$ are nonempty and consists solely of constant functions.
  Furthermore, $\alpha(u) \cap \omega(u) = \emptyset$, unless $u$ is itself constant.
\end{mythm}
\begin{myproof}
  Let $u$ be a bounded solution of \eqref{eq:grad_like} and consider the orbit
  \[
  \calO := \bigcup_{\tau\in\R} \{ \phi_\tau u \}.
  \]
  By assumption, the set $\calO$ is bounded in $C^{R-1}_b(\R,\R^d)$.
  Using the hypothesis \ref{hyp:G}, an elliptic bootstrapping argument allows us to conclude that $\calO$ is bounded in $C^{R+1}_b(\R,\R^d)$.
  Therefore, by the Arzel\`a--Ascoli theorem, $\calO$ is precompact in $C^R_{\text{loc}}(\R,\R^d)$.
  
  Now note that, since $K \in L^1(\R,\calM(\R^d))$, the map
  \[
  u \mapsto K*S(u) \; : \; E \cap C^R_{\text{loc}}(\R,\R^d) \mapsto C^R_{\text{loc}}(\R,\R^d)
  \]
  is continuous, for any bounded subset $E \subset C^R_b(\R,\R^d)$.
  Consequently, the map
  \[
  u \mapsto - \partial_x \nabla_v E(u,u') + \nabla_u E(u,u') + DS(u)^T K*S(u)) - c \Gamma(u) u'
  \]
  is continuous from $E \cap C^R_{\text{loc}}(\R,\R^d)$ into $C^0_{\text{loc}}(\R,\R^d)$.
  Hence, in particular, the accumulation points of $\calO$ in $C^R_{\text{loc}}(\R,\R^d)$
  are again $C^R_b(\R,\R^d)$-bounded solutions of \eqref{eq:grad_like}.
  We denote the closure of $\calO$ in $C^R_{\text{loc}}(\R,\R^d)$ by $\cl{\calO}$.
  
  We are now in a setup where we have a continuous flow $(\phi_\tau)_{\tau \in \R}$ on a compact metric space $\cl \calO$,
  where the flow possesses a weak Lyapunov function $\calL$.
  Compactness of $\cl \calO$ ensures that $\omega(u)$ is nonempty.
  Since $\calL$ is monotone along the flow, it follows that it is constant on $\omega(u)$.
  The $\omega$-limit set is flow invariant, that is, if $v \in \omega(u)$ then $\phi_\tau v \in \omega(u)$ for all $\tau \in \R$.
  Since $\calL$ is constant on $\omega(u)$, we have
  \[
  \ddxe{\tau} \calL(\phi_\tau v) = - c \Gamma(v) v'(\tau) \cdot v'(\tau) = 0, \qquad \text{for all} \quad \tau \in \R.
  \]
  This tells us that $\omega(u)$ consist solely of constant functions.
  The same argument applies to $\alpha(u)$.
  Clearly, if $u$ is not constant, then
  \[
  \calL(\omega(u)) = \lim_{\tau\to +\infty} \calL(\phi_\tau u) < \lim_{\tau\to -\infty} \calL(\phi_\tau u) = \calL(\alpha(u)),
  \]
  hence in that case $\alpha(u)$ and $\omega(u)$ must be disjoint.
\end{myproof}

\subsection{Extension to the center manifold}\label{sec:glcmfd}
We now consider $c \neq 0$ as a parameter in \eqref{eq:grad_like}, and we expand a parameter dependent center manifold around $c = 0$.
We allow for the maps $K$, $E$, $S$ to depend on additional parameters $\mu \in (-\epsilon',\epsilon')^s$ as well.
Define the principal symbol $\Sigma_p : \C \times (-\epsilon',\epsilon')^{s+1} \to \calM(\R^d)$ of \eqref{eq:grad_like} at $u=0$ by
\begin{equation}
  \label{eq:principal_symbol_grad}
  \Sigma_p(\nu ; c , \mu) := - \nabla_v^2 E_\mu(0,0) \nu^2 - c \Gamma(0) \nu + \nabla_u^2 E_\mu(0,0) + D^2 S_\mu(0)^T \kappa_0 S_\mu(0),
\end{equation}
where $\kappa_0 := \int_\R K_0(r) \d r$.
Likewise, the full symbol $\Sigma : \C \times (-\epsilon',\epsilon')^{s+1} \to \calM(\R^d)$ of \eqref{eq:grad_like} at $u=0$ is defined by
\begin{equation}
  \label{eq:full_symbol_grad}
  \Sigma(\nu ; c , \mu) := \Sigma_p(\nu ; c,\mu) + DS_\mu(0)^T \widehat K_\mu(\nu) DS_\mu(0).
\end{equation}
We then require the following.
\begin{myhyp}[C2$'$]
\namedlabel{hyp:C2_grad}{C2$'$}
  We demand that $\det \Sigma_p(\nu ; c,\mu) \neq 0$ for all $\nu \in i \R$, $(c,\mu) \in (-\epsilon',\epsilon')^{s+1}$.
  Furthermore, the map
  \[
  (\lambda,c,\mu) \mapsto  |\lambda| \| \Sigma_p(i \lambda ; c,\mu)^{-1} \|, \qquad \text{with} \quad \lambda \in \R, \quad (c,\mu) \in (-\epsilon',\epsilon')^{s+1},
  \]
  has to be uniformly bounded.
  If for some $\mu \in (-\epsilon',\epsilon')^{s+1}$ the map $E_\mu(u,v)$ is not constant in $v$, we require that the map
  \[
  (\lambda,c,\mu) \mapsto  |\lambda|^2 \| \Sigma_p(i \lambda ; c,\mu)^{-1} \|, \qquad \text{with} \quad \lambda \in \R, \quad (c,\mu) \in (-\epsilon',\epsilon')^{s+1},
  \]
  is uniformly bounded.
\end{myhyp}
Note that this is essentially a localized version of hypothesis \ref{hyp:G}.
Thus, under hypothesis \ref{hyp:C2_grad}, we find that small amplitude solutions of \eqref{eq:grad_like} exhibit gradient-like behavior.

In addition to hypothesis \ref{hyp:C2_grad}, we require that hypothesis \ref{hyp:Cmu} is satisfied.
We then obtain an unfolding of a parameter dependent center manifold $\calM_0(c,\mu)$, for $|c| , |\mu| < \epsilon'$.
Analogous to theorem \ref{thm:symplectic_center_manifold}, the Lyapunov functions $\calL_\mu$ can be restricted $C^k$ smoothly to $\calM_0(c,\mu)$.



\section{Applications}
\label{sec:applications}
We study the four examples \ref{ex:Allen-Cahn}--\ref{ex:triggered_GL} with the methods developed in the previous three chapters. 

\subsection{Periodic stationary patterns in nonlocal Allen-Cahn equations}
\label{sec:application_AC}

Recall from Example \ref{ex:Allen-Cahn} the stationary nonlocal Allen-Cahn equation
\begin{equation}
  \label{eq:AC_stationary}
  -\kappa_0 u + K*u - \nabla F(u) = 0,
\end{equation}
where $u : \R \to \R^d$, which possesses  a variational structure, with Lagrangian
\[
L(u,n) = E(u) + \frac 1 2 u \cdot n, \qquad E(u) = - \frac 1 2 u \cdot \kappa_0 u - F(u).
\]

\paragraph{Hamiltonian dynamics.}
Formally,  \eqref{eq:AC_stationary} exhibits a Hamiltonian structure, given by the Hamiltonian
\[
    \calH(u) = \rst{ \bigg( F(u) + \frac 1 2 u \cdot \kappa_0 u - \frac 1 2 u \cdot K * u \bigg) }{ x=0 } + \frac 1 2 \iint_Q u(x) \cdot K(x-y) u'(y) - u'(x) \cdot K(x-y) u(y) \d x \d y,
\]
where $Q=(-\infty,0) \times (0,\infty)$, and the presymplectic structure is given by
\[
\omega(v,w) = \iint_Q v(x) \cdot K(x-y) w(y) - w(x) \cdot K(x-y) v(y) \d x \d y.
\]

\paragraph{Center manifold reduction.}
Assume that the kernel $K$ is exponentially localized,
and that the nonlinearity $F\in C^{k+3}$, $k \geq 2$.
Moreover, suppose the following.
\begin{myhyp}[Nondegenerate zero solution]
We have
\[
\nabla F(0) = 0, \qquad \text{and} \quad \det\big( \kappa_0 + \nabla^2 F(0) \big) \neq 0.
\]
\end{myhyp}
The principal symbol $\Sigma_p$  \eqref{eq:principal_symbol} and  full symbol $\Sigma$ \eqref{eq:full_symbol} are
\[
\Sigma_p(\nu) = - \kappa_0 - \nabla^2 F(0),\qquad
\Sigma(\nu) = \widehat K(\nu) - \kappa_0 - \nabla^2 F(0).
\]
Let $i \ell_1,\dots,i \ell_N,-i \ell_1,\dots,-i \ell_N$ denote the roots, repeated according to multiplicity, of the characteristic equation
\begin{equation}
  \label{eq:AC_characteristic}
  \det \Sigma(\nu) = 0, \qquad \text{where} \quad \nu \in i \R.
\end{equation}
We then impose the following restriction.
\begin{myhyp}[Nonresonant eigenvalues]
  Assume the roots of the characteristic equation \eqref{eq:AC_characteristic} are nonresonant, that is,
  there exists no $0 \neq \n \in \Z^N$ such that $(\ell_1,\dots,\ell_N) \cdot \n = 0$.
\end{myhyp}
Under these conditions, hypothesis \ref{hyp:Comega} from section \ref{sec:Hamiltonian} is satisfied. We note that $\dim \ker \Sigma(\nu_i) = 1$ for all $1 \leq i \leq N$.
Choose $0 \neq \e_i \in \ker \Sigma(\nu_i)$, for each $1 \leq i \leq N$.
A basis for $\calE_0 := \ker \calT$ is now given by
\[
\psi_i(x) = \sin( \ell_i x ) \e_i, \qquad \psi_i^*(x) = \cos( \ell_i x ) \e_i, \qquad 1 \leq i \leq N.
\]
We choose a suitable projection $\calQ$ onto $\calE_0 \subset W^{1,2}_{-\eta}(\R,\R^d)$ as in the center manifold theorem, 
for example, an orthogonal projection in a weighted $L^2$ metric.
We obtain a center manifold $\calM_0 \subset W^{1,2}_{-\delta}(\R,\R^d)$ as the graph of a $C^k$ map
$\Psi : \calE_0 \to \ker \calQ$. On the center subspace $\calE_0$ we obtain a Hamiltonian $H$ and local symplectic form  $\scrw$ 
\[
H(u_0) := \calH\big( u_0 + \Psi(u_0) \big), \qquad 
\scrw_{u_0}(v,w) := \omega\big( v + D\Psi(u_0) v , w + D\Psi(u_0) w \big), \qquad \text{for} \quad u_0, v, w \in \calE_0.
\]
Locally, we obtain the induced Hamiltonian flow $\chi_H^\tau$ on $\calE_0$, defined by
\[
\partial_\tau \chi_H^\tau = V_H \circ \chi_H^\tau, \qquad \text{where} \quad - D H(u_0) = \scrw_{u_0}( V_H(u_0) , \cdot ).
\]
The flow of $\chi_H^\tau$ is conjugated, via the map $\id + \Psi$, to the flow of $\phi_\tau$ on $\calM_0$.

\paragraph{Existence of periodic solutions.}
We note that
\[
V_H(u_0) = \rst{ \partial_\tau \calQ\big( \phi_\tau u_0 + \phi_\tau \Psi(u_0) \big) }{\tau=0} = \calQ\big( u'_0 + \rst{\partial_\tau \phi_\tau \Psi(u_0)}{\tau=0} \big),
\]
where we use that the flow $\psi_\tau$ on $\calM_0$ is differentiable.
Consequently,
\[
D V_H(0) w = \calQ\big( w' + \rst{ \partial_\tau \phi_\tau D\Psi(0) w }{\tau=0} \big) = w', \qquad \text{for any} \quad w \in \calE_0,
\]
where we used that $D \Psi(0) = 0$.
Hence
\[
D V_H(0) = \partial_x : \calE_0 \to \calE_0.
\]
We then see that the eigenvalues and their algebraic multiplicity of the linearization ${D V_H(0)}$ coincide with the roots of \eqref{eq:AC_characteristic}.
The non-resonance condition then allows us to apply the Lyapunov center theorem, see \cite{abraham1978foundations}.
The theorem states that for each $1 \leq i \leq N$, there exists a family $u^{i,s}_0 \in \calE_0$, $s \in (0,\epsilon)$,
of small periodic solutions of the Hamiltonian flow $\chi_H^\tau$, such that
\[
\lim_{s\to 0} \sup_{\tau\in\R} \|\chi_H^\tau(u^{i,s}_0)\|_{\calE_0} = 0, \qquad \qquad \lim_{s\to 0} \period( u^{i,s}_0 ) = \frac{2\pi}{\ell_i}.
\]
In light of inequality \eqref{eq:Morrey_exponential}, the corresponding periodic orbits
\[
u^{i,s} := u^{i,s}_0 + \Psi( u^{i,s}_0 )
\]
on $\calM_0$ will be of small amplitude, hence solutions of the unmodified equation \eqref{eq:AC_stationary}. 

In a similar spirit to the result from section \ref{subsec:Whitham}, when a pair of eigenvalues coalesce, 
under suitable Krein signature and nondegeneracy conditions, 
the birth of a homoclinic orbit through a Hamiltonian-Hopf bifurcation could be detected.
It is conceivable that, under suitable smoothness and twist assumptions, existence of invariant tori could be established in a KAM type theorem.


\subsection{Traveling fronts in neural field equations}
\label{subsec:NFE}

Recall the equation for traveling waves in the neural field equation  from Example \ref{ex:NFE},
\begin{equation}
  \label{eq:NFE_TW}
  -c u' = - u + K*S_\mu(u).
\end{equation}
We consider here the scalar case, $u : \R \to \R$, $S_\mu : \R\to\R$, and $c \neq 0$.
The convolution kernel $K$ is assumed to be symmetric and exponentially localized, that is, $K(r) = K(-r)$ and $K \in W^{1,1}_{\eta_0}(\R,\R)$ for a sufficiently small $\eta_0 > 0$.
Furthermore, we assume that
\[
\kappa_0 := \int_\R K(r) \d r \neq 0.
\]
The parameter dependent nonlinearity should be of the following form.
\begin{myhyp}[Supercritical pitchfork bifurcation]
The nonlinearity $S_\mu$ is $C^{k+3}$ smooth, $k \geq 2$, and of the form
\[
S_\mu(u) = u \big( \mu + \kappa_0^{-1} - \alpha u^2 \big) + \calO\left( u^2 \big( |\mu| + |u| \big)^2 \right), \qquad \text{as} \quad (u,\mu) \to (0,0),
\]
with $\alpha > 0$.
\end{myhyp}
With this assumption, the constant solutions of \eqref{eq:NFE_TW} undergo a supercritical pitchfork bifurcation as $\mu$ passes through $0$.

\paragraph{Gradient-like behavior.}
Define smooth cutoff functions $\chi_\epsilon : \R \to \R$ by setting
\[
\chi_1(x) :=
\begin{cases}
  1 \quad &\text{if } |x| \leq 1, \\
  0 \quad &\text{if } |x| \geq 2,
\end{cases}
\]
and
\[
\chi_\epsilon(x) := \chi_1(x/\epsilon).
\]
Then define a modified nonlinearity $s_\mu$ by
\begin{equation}
  \label{eq:NFE_modified_S}
  s_\mu(u) := \chi_\epsilon(u) S_\mu(u)  -  \epsilon \big( 1-\chi_\epsilon(u) \big) u.
\end{equation}
Choosing $\epsilon > 0$ sufficiently small, we have, for some $\epsilon' > 0$,
\[
s_\mu'(u) < 0, \qquad \text{for all} \quad u \in \R,\; \mu \in (-\epsilon',\epsilon').
\]
Then, for small amplitude solutions, \eqref{eq:NFE_TW} is equivalent to
\begin{equation}
  \label{eq:NFE_TW_modified}
  -c s_\mu'(u) u' = - s_\mu'(u) u + s_\mu'(u) K* s_\mu(u).
\end{equation}
This is of dissipative form as in equation \eqref{eq:grad_like}, with
\begin{align*}
  L_\mu(u,n) &= E_\mu(u) + \frac 1 2 s_\mu(u)n, \qquad \text{where} \quad E_\mu(u) = - \int_0^u s_\mu'(r) r \d r, \\
  \Gamma(u) &= - s_\mu'(u).
\end{align*}
Hence, by lemma \ref{lemma:Lyapunov} we obtain a Lyapunov function
\begin{multline}
  \label{eq:NFE_Lyapunov}
    \calL_\mu(u) = - \rst{ \bigg( E_\mu(u) + \frac 1 2 s_\mu(u) K * s_\mu(u) \bigg) }{x=0} \\
  + \frac 1 2 \iint_Q K(x-y) \bigg( s_\mu(u(x)) s_\mu'(u(y)) u'(y) - s_\mu(u(y)) s_\mu'(u(x)) u'(x) \bigg) \d x \d y.
\end{multline}
In light of theorem \ref{thm:grad_like}, solutions of the modified equation \eqref{eq:NFE_TW_modified} have gradient-like behavior.
Since the modified equation \eqref{eq:NFE_TW_modified} locally coincides with the original equation \eqref{eq:NFE_TW},
we conclude that for small values of $\mu$, small amplitude solutions of \eqref{eq:NFE_TW} (with $c \neq 0$) exhibit gradient-like behavior.

\paragraph{Center manifold reduction.}
The principal and full symbol from \eqref{eq:principal_symbol_grad} and  \eqref{eq:full_symbol_grad}, respectively, are
\[
\Sigma_p(\nu ; c,\mu) = \big( \mu + \kappa_0^{-1} \big) \big( c \nu - 1 \big),\qquad \Sigma(\nu ; c,\mu) = \big( \mu + \kappa_0^{-1} \big) \big( c \nu - 1 \big) + \frac 1 2 \big( \mu + \kappa_0^{-1} \big)^2 \widehat K(\nu).
\]
We readily see that hypothesis \ref{hyp:Cmu} is satisfied.

Note that $\nu = 0$ is always a double root of the characteristic equation $\Sigma(\nu ; 0,0) = 0$.
Now impose the following restriction upon $K$.
\begin{myhyp}[Unfolding of the kernel]
Assume that
\[
\Sigma(\nu;c,\mu) = ( \nu^2 + c \nu - \mu ) d(\nu ; c , \mu),
\]
where $d(\nu ; c , \mu) \neq 0$ for all $\nu \in i \R$ and $|c| , |\mu| < \epsilon'$, for some $\epsilon' > 0$.
\end{myhyp}

A basis for $\calE_0 := \ker \calT_{0,0}$ is now given by
\[
\e_0(x) = 1, \qquad \e_1(x) = x.
\]
We define the projection $\calQ$ onto $\calE_0$ by
\[
\calQ u := u(0) \e_0 + u'(0) \e_1.
\]
This projection extends to $W^{2,2}_{-\eta}(\R,\R^d)$, independent of $\eta > 0$.
We will unfold a parameter dependent center manifold $\calM_0(c,\nu) \subset W^{2,2}_{-\delta}(\R,\R^d)$
as the graph of a map
\[
\Psi : \calE_0 \times (-\epsilon',\epsilon')^2 \to \ker \calQ.
\]

We henceforth identify $\calE_0$ with $\R^2$, via the map
\[
(A,B) \mapsto u_0 = A \e_0 + B \e_1.
\]
Through pullback along the map $\id + \Psi$, 
the shift map $\phi_\tau$ on $\calM_0(c,\mu)$ induces a parameter dependent smooth flow $\chi^\tau_{c,\mu}$ on $\R^2$.
The function $\Lambda_{c,\mu} : \R^2 \to \R$ defined by
\[
\Lambda_{c,\mu}(A,B) := \calL_\mu\big( u_0 + \Psi(u_0,c,\mu) \big), \qquad u_0 = A \e_0 + B \e_1,
\]
is a Lyapunov function for the flow $\chi^\tau_{c,\mu}$, whenever $c \neq 0$.

\paragraph{Existence of traveling fronts.}
We will make use of Conley index theory, see for example \cite{MischaikMrozek} for an overview.
First, to compute the index we will employ the Conley continuation theorem along the parameter $\mu$.
This relies on the existence of an isolating neighborhood, uniform in the parameter $\mu$.
To obtain such a neighborhood, we compute
\[
    \Lambda_{c,\mu}(A,B) = - \frac 1 4 \alpha A^4 - \frac 1 4 \frac{\kappa_2}{\kappa_0^2} B^2  + \scro \left( \big( A^2 + |B| \big)^2 \right) 
  + \calO\left( \big( |c| + |\mu| \big) \big( |A| + |B| \big)^2 \right)
\]
as $(A,B,c,\mu) \to (0,0,0,0)$; see Appendix \ref{appendix:NFE} for details.
We then impose the following restriction upon $K$.
\begin{myhyp}[Positive second moment]
  We require $\kappa_2 > 0$.
\end{myhyp}
Using the expansion of $\Lambda_{c,\mu}$, we then find that
there exists a neighborhood $0 \in V \subset \R^2$ such that for any $\mu \in (-\epsilon',\epsilon')$ the flow $\chi_{c,\mu}^\tau$ is transverse to $\bdy V$.
We again refer to Appendix \ref{appendix:NFE} for details of this construction.
We can now employ the Conley continuation theorem, which states that
\[
\HC_k\big( V, \chi^\tau_{c,\mu_1} \big) \iso \HC_k\big( V , \chi^\tau_{c,\mu_2} \big), \qquad \text{for any} \quad \mu_1 , \mu_2 \in (-\epsilon',\epsilon').
\]

For $\mu < 0$, the set $V$ contains precisely one stationary point of the flow $\chi^\tau_{c,\mu}$, namely, $(A,B) = (0,0)$.
It then follows from the gradient-like behavior that
\[
\HC_k\big( V , \chi^\tau_{c,\mu} \big) \iso \HC_k\big(\{(0,0)\} , \chi^\tau_{c,\mu} \big), \qquad \text{for} \quad \mu < 0. 
\]
For small $\mu >0$, the set $V$ contains precisely three stationary points of the flow $\chi^\tau_{c,\mu}$, say, $(A,B)=(0,0)$ and $(A,B) = ( A^\pm_\mu , 0 )$.
Suppose however that it does not contain a heteroclinic solution between $(0,0)$ and one of the points $(A_\mu^-,0)$, $(A_\mu^+,0)$.
Then, in light of the gradient-like behavior of \eqref{eq:NFE_TW_modified}, we find that the invariant set contained in $V$ is
\[
\Inv(V) = \{ (0,0) \} \cup I, \qquad \text{where} \quad \{ (0,0) \} \cap I = \emptyset.
\]
Here, the set $I$ could consist of only $(A_\mu^-,0)$ and $(A_\mu^+,0)$, or it could contain a heteroclinic between them.
By additivity of the Conley index, for $\mu > 0$ we would then have
\[
\HC_k\big( V, \chi^\tau_{c,\mu} \big) \iso
\HC_k\big( \{(0,0)\} , \chi^\tau_{c,\mu} \big) \oplus \HC_k\big( I , \chi^\tau_{c,\mu} \big).
\]
Combining these results, we would find that
\begin{equation}
  \label{eq:indices_contradiction}
  \HC_k\big( \{(0,0)\}, \chi^\tau_{c,-\mu} \big) \iso \HC_k\big( \{(0,0)\} , \chi^\tau_{c,\mu} \big) \oplus \HC_k\big( I , \chi^\tau_{c,\mu} \big).
\end{equation}
Now let us without loss of generality assume that $c > 0$.
Our hypothesis on the destabilization of $0$ then implies that $(0,0)$ is stable for $\mu < 0$, and a saddle for $\mu > 0$.
But then, for small $\mu > 0$ we have
\[
\HC_k\big( \{(0,0)\}, \chi^\tau_{c,-\mu} \big) \iso
\begin{cases}
  \Z_2 \quad & \text{if } k = 0, \\
  0 \quad & \text{otherwise},
\end{cases}
\qquad
\HC_k\big( \{(0,0)\}, \chi^\tau_{c,\mu} \big) \iso
\begin{cases}
  \Z_2 \quad & \text{if } k = 1, \\
  0 \quad& \text{otherwise}.
\end{cases}
\]
This contradicts \eqref{eq:indices_contradiction}.
We may therefore conclude,
the existence of at least one heteroclinic point $(A,B)=(A_h,B_h)$  for any given small positive wave speed $c > 0$ and small $\mu > 0$, such that the heteroclinic orbit $\bigcup_\tau \chi_{c,\mu}^\tau(A_h,B_h)$ is contained in $V$,
and connects $(0,0)$ with one of the points $(A_\mu^-,0)$, $(A_\mu^+,0)$.
The definition of the trapping region $V$ ensures that the heteroclinic
\[
u_h = A_h \e_0 + B_h \e_1 + \Psi(A_h,B_h,c,\mu)
\]
satisfies
\[
\sup_{\tau\in\R} \| \phi_\tau u_h \|_{W^{1,2}_{-\delta}(\R,\R^d)} < \frac{\epsilon}{C_{0,2,-\eta}}.
\]
Hence by inequality \eqref{eq:Morrey_exponential} we have $\| u_h \|_{C^0_b(\R,\R^d)} < \epsilon$, so that $u_h$ is a solution of the unmodified equation \eqref{eq:NFE_TW}.

If, in addition, we assume that $s_\mu$ is an odd function for each $\mu \in (-\epsilon',\epsilon')$, we have the $\Z_2$-symmetry $u(x) \mapsto -u(x)$ in \eqref{eq:NFE_TW_modified}.
Then $A^-_\mu = - A^+_\mu$, hence $\Lambda_{c,\mu}(A^-_\mu) = \Lambda_{c,\mu}(A^+_\mu)$, hence there cannot exist a heteroclinic orbit connecting $(A^-_\mu,0)$ and $(A^+_\mu,0)$.
Furthermore, by symmetry considerations, there must then exist at least two heteroclinic solutions $u^\pm_h$ of \eqref{eq:NFE_TW} for small $c \neq 0$, $\mu > 0$,
one connecting $(0,0)$ with $(A_\mu^-,0)$, and one connecting $(0,0)$ with $(A_\mu^+,0)$.

Similar results can be obtained for the vector-valued neural field equation, adapting the assumptions:
\begin{enumerate}
\item For any $u$, the linearization $DS(u)$ is positive definite.
\item There exists a function $E : \R^d \to \R$ such that $\nabla E(u) = - DS(u) u$.
\item The full symbol $\Sigma(u;c,\mu)$ satisfies Hypothesis \ref{hyp:C4}.
\end{enumerate}


\subsection{Solitons and periodic waves in Whitham type equations}
\label{subsec:Whitham}
We recall from Example \ref{ex:Whitham} that traveling waves $u(t,x) = u(\xi)$, $\xi = x-ct$  in Whitham type equation
\[
u_t + 2 \alpha u u_x + \int K(x-y) u_x(t,y) \d y = 0
\]
satisfy
\[
\bigg( \alpha u^2 - c u + K * u \bigg)' = 0.
\]
Integrating and assuming $u\to 0$ at infinity, one obtains
\begin{equation}
  \label{eq:Whitham_classic}
  \alpha u^2 - c u + K * u = 0.
\end{equation}
The question of existence of homoclinic and periodic solutions to \eqref{eq:Whitham_classic} is considered in for example \cite{ehrnstrom2009traveling, ehrnstrom2012existence} using a variational approach.
We present here an alternative proof for the existence of small amplitude homoclinics and periodic solutions, using the Hamiltonian formalism for nonlocal equations. A classical choice for the kernel is (see \cite{whitham2011linear})
\begin{equation}
  \label{eq:Whitham_kernel}
  K(x) = \int_\R e^{i x \xi} \sqrt{ \frac{ \tanh(\xi) }{ \xi } } \d \xi.
\end{equation}
Unfortunately, the formulation of the results in \cite{faye2013fredholm, faye2016center} does not allow for the kind of mild singularity at the origin that this kernel exhibits. Although, upon inspection, these assumptions there could in fact be weakened to include such singularities, we will restrict ourselves here to Whitham type equations where the kernel satisfies hypothesis \ref{hyp:C3}. The problem \eqref{eq:Whitham_classic} can be written in the slightly more general form 
\begin{equation}
  \label{eq:Whitham_TW}
  K * u + F'_{c}(u) = 0,
\end{equation}
for some parameter dependent potential $F_{c} : \R \to \R$. Motivated by \eqref{eq:Whitham_classic} and \eqref{eq:Whitham_kernel}, we impose the following restrictions upon $K$ and $F_c$.
\begin{myhyp}[Critical wave speed]
There exists a $c_* \neq 0$ so that the following holds.
\begin{enumerate}
\item 
  The kernel $K$ satisfies
  \[
  \widehat K(\nu) - c = ( -\nu^2 + c-c_* ) d(\nu ; c), \qquad \text{with} \quad \widehat K(\nu) = \int_\R e^{-\nu x} K(x) \d x,
  \]
  where $d(\nu ; c) \neq 0$ for all $\nu \in i \R$ and $|c-c_*| < \epsilon'$, for some $\epsilon' > 0$.

  \item
The nonlinearity $F_c$ is $C^{k+3}$ smooth, $k \geq 2$, and of the form
\[
F_c(u) = - \frac 1 2 c u^2 + \frac{ 1  }{ 3 } \alpha u^3 + \calO\left( u^2 \big( |u| + |c-c_*| \big)^2 \right), \qquad \text{as} \quad (u,c) \to (0,c_*),
\]
for some $\alpha \neq 1/3$.
\end{enumerate}
\end{myhyp}
A quick calculation shows that the speed $c_*$ is the group velocity of linear waves near $u=0$.

\paragraph{Hamiltonian dynamics.}
Equation \eqref{eq:Whitham_TW} has a variational structure, with Lagrangian, Hamiltonian, and presymplectic structure,
\begin{equation}
  \label{eq:Whitham_Hamiltonian}
  \begin{split}
   L_c(u,n) &= F_c(u) + \frac 1 2 u n,\\
  \calH_c(u) &= - \rst{ \bigg(F_c(u)  + \frac 1 2 u K * u \bigg) }{ x=0 } 
 + \frac 1 2 \iint_Q K(x-y)  \bigg( u(x) u'(y) - u'(x) u(y) \bigg) \d x \d y,\\
  \omega(v,w)& = \iint_Q K(x-y) \bigg( v(x) w(y) - w(x) v(y) \bigg) \d x \d y, \qquad Q=(-\infty,0)\times (0,\infty).
  \end{split}
\end{equation}

\paragraph{Center manifold reduction.}
Assume that the kernel $K$ is exponentially localized.
The principal and full  symbol $\Sigma_p$ and $\Sigma$ from \eqref{eq:principal_symbol} and \eqref{eq:full_symbol}, are 
\[
\Sigma_p(\nu ; c) = - c,
\qquad \qquad \Sigma(\nu ; c) = \widehat K(\nu) - c,
\]
thus satisfying Hypothesis \ref{hyp:Cmu} near $c = c_*$. A basis for $\calE_0 := \ker \calT_0$  and a projection $\calQ$ onto $\calE_0$ are  
\[
\e_0(x) = 1, \qquad \e_1(x) = x,\qquad 
\calQ u := u(0) \e_0 + u'(0) \e_1.
\]
The projection extends to $W^{2,2}_{-\eta}(\R,\R^d)$, for all $\eta > 0$.
We will unfold a parameter dependent center manifold $\calM_0(c-c_*) \subset W^{2,2}_{-\delta}(\R,\R^d)$
as the graph of a map
\[
\Psi : \calE_0 \times (-\epsilon',\epsilon') \to \ker \calQ.
\]
and choose coordinates in $\calE_0$  via 
\[
(A,B) \mapsto u_0 = A \e_0 + B \e_1.
\]
Through pullback along the map $\id + \Psi$, 
the shift map $\phi_\tau$ on $\calM_0(c - c_*)$ induces a parameter dependent smooth flow $\chi^\tau_{c}$ on $\R^2$.
The function $H_{c} : \R^2 \to \R$ defined by
\[
H_{c}(A,B) := \calH_{c}\big( u_0 + \Psi(u_0,c-c_*) \big), \qquad u_0 = A \e_0 + B \e_1,
\]
is the Hamiltonian for the flow $\chi^\tau_{c}$.
The corresponding symplectic form $\widetilde \omega_{(A,B),c}$ on $\R^2$ is defined via the pullback relation
\[
\widetilde \omega_{(A,B),c}( v , w ) := \omega\big( v + D_u \Psi(u_0,c-c_*) v , w + D_u \Psi(u_0,c-c_*) w \big), \qquad \text{where} \quad v , w \in \R^2.
\]

We note that a parameter independent center manifold can be expanded around $\mu = \mu_*$ for each $c \geq 0$.
Hence the center manifold exists for all $c \in (- \epsilon' , \infty)$ and $\mu$ close to $\mu_*$,
and the flow of small amplitude $u \in \calM_0(c,\mu-\mu_*)$ is Hamiltonian, for all $c \geq c_*$.

\paragraph{Asymptotic expansion of $H_c$.}
By expanding the center manifold and the Hamiltonian, one computes
\begin{multline*}
  H_c(A,B) = \frac{1-3\alpha}{6} A^3 + \frac 1 2 (c-c_*) A^2 - \frac 1 4 \kappa_2 B^2 \\ +\calO\left( A^4 + |B|^3 + |A| B^2+ |A|^3 |B|+|c-c_*|( |A|^3 + |AB| + B^2 )+ |c-c_*|^2 A^2  \right)
\end{multline*}
as $(A,B,c) \to (0,0,c_*)$; see Appendix \ref{appendix:Whitham} for a detailed computation.
We define the truncated Hamiltonian
\[
\widetilde H_c(A,B) = \frac{1-3\alpha}{6} A^3 + \frac 1 2 (c-c_*) A^2 - \frac 1 4 \kappa_2 B^2.
\]
In the following, we discuss the reduced flow for the case  $\kappa_2 > 0$, as is the case when $K$ is of the form \eqref{eq:Whitham_kernel}. The case $\kappa_2 < 0$  yields similar results.

\paragraph{Existence of solitons for $c > c_*$.}
For $c > c_*$, the geometry of the level set $\set{ (A,B) }{ \widetilde H_c(A,B) = \widetilde H_c(0,0) }$ dictates the existence of a homoclinic point $(A_h,0)$
in the modified flow $\widetilde \chi_c^\tau$ on $\R^2$, defined through
\[
\partial_\tau \widetilde \chi_c^\tau = V_{\widetilde H_c} \circ \widetilde \chi_c^\tau, \qquad \text{where} \quad - D \widetilde H_c(u) = \widetilde \omega_{(A,B),c}(V_{\widetilde H_c}(A,B) , \cdot ).
\]
Define the homoclinic orbit
\[
\widetilde \gamma_h^c := \{0\} \cup \bigcup_{\tau\in\R} \big\{ \widetilde \chi_c^\tau(A_h,0) \big\}.
\]
One readily estimates
\[
\dist( 0 , \widetilde \gamma_h^c ) = \calO\left( |c-c_*| \right) \qquad \text{as} \quad c \to c_*.
\]
Consequently, the unperturbed Hamiltonian vector field $V_{H_c}$ is a $C^1$-small perturbation of $V_{\widetilde H_c}$.
Furthermore, note the homoclinic intersects the plane $B=0$ transversely.
By the reversible symmetry $(A,B) \mapsto (A,-B)$, inherited from the symmetry $u(x) \mapsto u(-x)$, the homoclinic will persist in the unperturbed system, for small $|c-c_*|$.
The smallness of the orbit $\widetilde \gamma_h^c$, combined with inequality \eqref{eq:Morrey_exponential}, implies the existence of a homoclinic $u_h^c$ to $0$ in the full system \eqref{eq:Whitham_TW}.

\paragraph{Periodic waves.}
In the regime $c > c_*$, it follows that the homoclinic orbit is shadowed by periodic orbits.
More precisely, using \cite{vanderbauwhede1992homoclinic}, we see that there exists a family $u^{c,s}$, $s \in (0,\epsilon)$ of periodic orbits, such that
\[
\lim_{s \to 0} \dist\big( u^{c,s} , \Gamma_h^c \big) = 0, \qquad \qquad \lim_{s\to 0} \period(u^{c,s} ) = \infty,
\]
where
\[
\Gamma_h^c := \bigcup_{\tau\in\R} \big\{ \phi_\tau u_h^c \big\}.
\]
For parameter values $c < c_*$, the linearization of the Hamiltonian flow at $0$ has eigenvalues $\pm i \sqrt{ |c - c_*| }$.
Arguing as in Section \ref{sec:application_AC}, one proves the existence of a family $u^{c,s} \in \calM_0$, $s \in (0,\epsilon)$, of periodic solutions of \eqref{eq:Whitham_TW}, such that
\[
\lim_{s\to 0} \|u^{c,s}\|_{C^0_b(\R,\R^d)} = 0, \qquad \qquad \lim_{s\to 0} \period( u^{c,s} ) = \frac{2\pi}{\sqrt{|c-c_*|}}.
\]

We remark here that similar results could have been obtained without the Hamiltonian structure,
but instead relying on local bifurcation analysis using reversibility, only.


\subsection{Wavenumber selection in NLS}

Recall from Example \ref{ex:triggered_GL} the triggered nonlinear, nonlocal Schr\"odinger equation
\begin{equation}
  \label{eq:triggered_GL}
  A'' + \left( \lambda(x) - \frac{c^2}{4} \right) A + DS(A) K* S(A) = 0.
\end{equation}
Here $S(u) = f( |u|^2 )$ and we assume $\lambda$ satisfies
\[
\lambda(x) =
\begin{cases}
  0 \quad &\text{if } x \leq -\ell, \\
  1 \quad &\text{if } x \geq \ell,
\end{cases}
\]
for some value of $\ell > 0$.
Using the canonical identification $\C \iso \R^2$, we rewrite \eqref{eq:triggered_GL} into a system of real equations,
\begin{equation}
  \label{eq:triggered_GL_R}
  A'' + \left( \lambda(x) - \frac{c^2}{4} \right) A + Ds(A)^T k* s(A) = 0.
\end{equation}
where $A : \R \to \R^2$,
\[
s(u) :=
\begin{pmatrix}
  f(|A|^2) \\
  f(|A|^2)
\end{pmatrix},
\qquad
k(r) :=
\begin{pmatrix}
  K(r) & 0 \\
  0 & K(r)
\end{pmatrix}.
\]
Then \eqref{eq:triggered_GL_R} has a variational structure, with Lagrangian
\[
L(x,u,v,n) = E(x,u,v) + \frac 1 2 s(u) \cdot n, \qquad \text{where} \quad E(x,u,v) = - \frac 1 2 |v|^2 + \frac 1 2 \left( \lambda(x) - \frac{c^2}{4} \right)  |u|^2.
\]

\paragraph{Symmetry and conserved quantity.}
The symmetry group of \eqref{eq:triggered_GL_R} is $G = \mathbf{O}(2) \times \{1\}$.
The Lie algebra $\frakg = \frako(2) \times \{0\}$ consists of matrices $\xi$ of the form
\begin{equation}
  \label{eq:GL_Lie_algebra}
  \xi =
\alpha
\begin{pmatrix}
  0 & 1 \\
  -1 & 0
\end{pmatrix},
\qquad \alpha \in \R,
\end{equation}
acting on functions $u : \R \to \R^2$ via left multiplication.
Hence the conserved quantity $\calC_\xi(u)$ as defined in theorem \ref{thm:Noether} is of the form
\[
\calC_\xi(u)(x) = \nabla_v E(x,u,u') \cdot u_\xi + B_\xi(\phi_x u).
\]
Since $u \mapsto S(u)$ is invariant under the $G$-action on $u$, we have 
\[
DS(u) u_\xi = \ddxe{\tau} S\bigg( \exp_e(\tau \xi) \bullet u \bigg) = 0.
\]
Hence $\calB_\xi(\phi_x u) = 0$.
Thus
\[
\calC_\xi(u) = - u'(x) \cdot \xi u(x).
\]
Note that this quantity is independent of $c$ and $\lambda$.

\paragraph{Restriction on wavenumbers.}
We consider the behavior of \eqref{eq:triggered_GL_R} with a constant potential $\mu$, that is,
\begin{equation}
  \label{eq:stationary_NLS}
  A'' + \mu A + DS(A) K* S(A) = 0,
\end{equation}
where we will have  $\mu = - c^2/4$ near $x=-\infty$ and $\mu=1-c^2/4$ near $x=\pm\infty$.
Plane waves ${A(x) = \exp_e( x \xi ) \bullet A_0}$, with $\xi \in \frako(2)$ and $A_0 \in \R^2\setminus\{0\}$, solve
\begin{equation}
  \label{eq:GL_planar_wave_condition} \|\xi\|^2 = \mu + 2 \kappa_0 f'{\left(|A_0|^2\right)} f{\left(|A_0|^2\right)},
\end{equation}
with $\kappa_0 = \int_\R K(r) \d r$.
%
Now suppose $A \in C^2_b(\R,\R^2)$ solves \eqref{eq:triggered_GL_R} and connects to plane waves $A_\pm$, that is,
$\alpha(A) = A_-$ and $\omega(A) = A_+$, with $\alpha$ and $\omega$ as in Definition \ref{def:limit_sets}.
Since $A_\pm$ are planar waves, we have $A_\pm' = \xi^\pm A_0^\pm$,
hence
\[
\calC_\xi(A_\pm)(x) = - | \xi^\pm A_0^\pm |^2.
\]
Since $\calC_\xi(u)$ must be constant, it must hold that
\begin{equation}
  \label{eq:conserved_angular_momentum}
  - | \xi^- A_0^- |^2 = \lim_{x\to-\infty} \calC_\xi(A)(x) = \lim_{x\to\infty} \calC_\xi(A)(x) = - | \xi^+ A_0^+ |^2.
\end{equation}

Now if, for example,
\[
2 \kappa_0 f'{\left(|a|^2\right)} f{\left(|a|^2\right)} < \frac{c^2}{4}, \qquad \text{for all} \quad a \in \R,
\]
it follows from \eqref{eq:GL_planar_wave_condition} that $A_- = 0$.
Hence, in light of \eqref{eq:conserved_angular_momentum}, either $\xi^+=0$ or $A_0^+ = 0$.
Likewise, if
\[
2 \kappa_0 f'{\left(|a|^2\right)} f{\left(|a|^2\right)} < \frac{c^2}{4} - 1, \qquad \text{for all} \quad a \in \R,
\]
we have $A_+ = 0$ and therefore either $\xi^- = 0$ or $A_0^+ = 0$.
In both cases, we find the \emph{selection} of the wavenumber $\xi=0$ for nontrivial plane waves compatible with the parameter step, based on conservation of the angular momentum associated with the gauge symmetry.




\section{Discussion}
We considered a class of nonlocal equations, where the nonlocal contribution is in the form of a symmetric convolution operator. With the aid of a nonlocal analogue of Green's formula, we derived conserved quantities out of symmetries of the equation. For translation symmetries $u(\cdot) \mapsto u(\cdot+\tau)$, the associated conserved quantity can be interpreted as the Hamiltonian of the nonlocal system, generating the ``dynamics'' through a symplectic formalism, at least formally, but rigorously on center manifolds. We also established gradient-like dynamics when dissipative terms, for instance advection terms from a traveling-wave ansatz, are added.  We comment here on possible generalizations and limitations to our approach.

\paragraph{More general nonlocal operators.}
The specific form of a bilinear form composed of pointwise evaluation and convolution is clearly  unnecessarily restrictive. Generalizations could include multilinear convolution, singular integral, or pseudo-differential operators. It is conceivable that the Noether theorem extends to such equations. Interesting, more general, examples of symmetries would then be nonlocal symmetries arising through scale invariance ${u(x) \mapsto \alpha(\tau) u(\beta(\tau) x)}$. Particular examples in this direction are nonlocal operators based on fractional powers of the Laplacian. Based on a local formulation, Hamiltonian identities have been derived in this case \cite{sire} and exploited in existence and uniqueness proofs \cite{silvestre}. The present techniques point towards generalizations beyond the particular, scale-invariant fractional Laplacian with its interpretation as a local operator in an extended half space.

\paragraph{Infinite-dimensional systems.}
We considered nonlocal equations where $u : \R \to \R^d$, only. Natural generalizations would extend to  $u : \R \to \calX$, where $\calX$ is a Hilbert or Banach space. Extending the Noether theorem does not appear to cause any difficulties and leads to generalized Hamiltonian identities as in \cite{gui2008hamiltonian}. We caution however that results based on center manifold reduction or Fredholm theory would require an extension of \cite{faye2013fredholm}.

\paragraph{Global topological invariants.}
The example in Section \ref{subsec:NFE} illustrates how the gradient-like behavior,
combined with robust topological machinery such as the Conley index,
can be used to prove the generic existence of heteroclinic orbits in dissipative systems.
We relied on a center manifold reduction, since Conley index theory is typically only applicable in finite dimensional settings. Extensions of these ideas to large amplitude solutions are conceivable, albeit much more technically involved. One would likely have to develop a new topological invariant, tailored towards dissipative nonlocal equations, which captures the dynamical properties.
The spectral flow theory from \cite{faye2013fredholm}, combined with the gradient-like behavior, gives rise to a setting in which one can analyze moduli spaces of connecting orbits.
It seems likely that this ultimately could lead to a global Floer homology theory, comparable to \cite{bakker2017floer, angenent1999superquadratic}.

\paragraph{Global Hamiltonian dynamics.}
It would be interesting and quite useful to find conditions on the kernel $K$ and the underlying function spaces,
such that $\omega$ is nondegenerate, without restricting to a center manifold.
This would put the nonlocal equations in the context of Hamiltonian PDEs, as considered for example in \cite{kuksin2000analysis}, possibly allowing for global Lyapunov center theorems and KAM theory. 
A first step in this direction could be an extension of center manifold theory \cite{faye2016center} to neighborhoods of periodic solutions,
comparable to \cite{dynOfNonlinearWaves},
where new questions arise, both in terms of regularity constraints necessary for the reduction procedure and in terms of the reduced dynamics.

%

%



\appendix
\section[Appendix]{Computations on the center manifold}
\label{sec:appendix}


\subsection{Neural field equations}
\label{appendix:NFE}

\subsubsection{Constructing a trapping region}
We construct a small neighborhood $V$ of $(0,0)$ for which the flow of $\chi^\tau_{c,\mu}$ is transverse to the boundary $\bdy V$,
for a fixed $c \neq 0$, but uniform in $\mu \in (-\epsilon',\epsilon')$, small. First, we formally expand 
\[
\Psi(A,B,c,\mu) = \sum_{i,j,k,l \geq 1} A^i B^j c^k \mu^l \psi_{i,j,k,l}, \qquad \text{where} \quad \psi_{i,j,k,l} \in \ker \calQ.
\]
The tangency $D_u \Psi(0,0,0,0) = 0$ implies $\psi_{1,0,0,0} = \psi_{0,1,0,0} = 0$.
For $c = \mu = 0$, equation \eqref{eq:NFE_TW_modified} is invariant under the map $u(x) \mapsto - u(x)$.
It follows that $\Psi(A,B,0,0) = - \Psi(-A,-B,0,0)$, such that $\psi_{2,0,0,0}=\psi_{0,2,0,0}=\psi_{1,1,0,0}=0$.

In order to compute the Taylor expansion  $\Lambda_{c,\mu}(A,B)$, we consider the boundary term in \eqref{eq:NFE_Lyapunov}, 
\begin{equation}
  \label{eq:NFE_boundary}
  \calB_\mu(u) = \frac 1 2 \iint_Q K(x-y) \bigg( s_\mu(u(x)) s_\mu'(u(y)) u'(y) - s(u(y)) s_\mu '(u(x)) u'(x) \bigg) \d x \d y.
\end{equation}
Formally, we let
\[
\calB_\mu(u) = \sum_{i,j,k,l \geq 1} A^i B^j c^k \mu^l b_{i,j,k,l}.
\]
First note that
\begin{equation}
  \label{eq:partial_expansion_s}
  s_\mu(u) s_\mu'(v) w = \kappa_0^{-2} u w - \alpha \kappa_0^{-1} u^3 w - 3 \alpha \kappa_0^{-1} u v^2 w + \calO\left( \big(|u|+|v|\big)^4 |w| + |\mu| |u| |w| \right)
\end{equation}
as $(u,v,w,\mu) \to (0,0,0,0)$.
In the following, we use the notations
\begin{equation}
  \label{eq:skew_pairing}
  P(u,v) := \frac 1 2 \iint_Q K(x-y) \bigg( u(x) v(y) - v(x) u(y) \bigg) \d x \d y,\qquad \kappa_m := \int_\R r^m K(r) \d r.
\end{equation}
We now compute terms up to order $\calO\left( \big( A^2 + |B| \big)^2 \right)$ in the expansion of $\calB_\mu$.
\begin{description}
\item[Terms of order $\calO( |A| + |B| )$.]
From \eqref{eq:partial_expansion_s} we see that the expansion of $\calB$ does not contain linear terms, 
hence $
b_{1,0,0,0} = b_{0,1,0,0} = 0.
$
\item[Terms of order $\calO(A^2)$.]
At order $\calO(A^2)$ we obtain
\[
b_{2,0,0,0} = \kappa_0^{-2} P( \e_0 , \e_0' ) = \kappa_0^{-2} P( \e_0 , 0 ) = 0.
\]
\item[Terms of order $\calO(B^2)$.]
At order $\calO(A^2)$ we find, using Lemma \ref{lemma:omega_covariance_computation},
\[
b_{0,2,0,0} = \kappa_0^{-2} P( \e_1 , \e_1' ) = \frac 1 2 \kappa_0^{-2} C_{1,0} \kappa_2 = - \frac 1 4 \frac{\kappa_2}{\kappa_0}.
\]
\item[Terms of order $\calO(|AB|)$.]
At order $\calO(|AB|)$ we find, using Lemma \ref{lemma:omega_covariance_computation} and $\kappa_1=0$ due to $K$ being even,
\[
b_{1,1,0,0} = \kappa_0^{-2} \bigg( P(\e_0,\e_1') + P(\e_1,\e_0') \bigg) = \kappa_0^{-2} P(\e_0,\e_1') = \frac 1 2 \kappa_0^{-2} C_{0,0} \kappa_1 = 0.
\]
\item[Terms of order $\calO(|A|^3)$.]
At order $\calO(|A|^3)$ we have,  using that $\psi_{2,0,0,0} = 0$,
\[
b_{3,0,0,0} = \kappa_0^{-2} P(\e_0 , \psi_{2,0,0,0}' ) = \kappa_0^{-2} P(\e_0,0) = 0.
\]
\item[Terms of order $\calO(A^2|B|)$.]
At order $\calO(A^2|B|)$ we find, using that$\psi_{2,0,0,0} = \psi_{1,1,0,0} = 0$,
\[
b_{2,1,0,0} = \kappa_0^{-2} \bigg( P(\e_1 , \psi_{2,0,0,0}' ) + P(\psi_{2,0,0,0},\e_1') + P( \e_0 , \psi_{1,1,0,0}' ) + P(\psi_{1,1,0,0},\e_0')\bigg) = 0.
\]
\item[Terms of order $\calO(A^4)$.]
Finally, at order $\calO(A^4)$ we obtain
\[
b_{4,0,0,0} = \kappa_0^{-2} P(\psi_{2,0,0,0},\psi_{2,0,0,0}') + 2 \alpha \kappa_0^{-1} P(\e_0^3,\e_0') - 6 \alpha \kappa_0^{-1} P( \e_0 , \e_0^2 \e_0' ) = 0.
\]
\end{description}
We conclude that
\[
\calB_\mu(u) = - \frac 1 4 \frac{\kappa_2}{\kappa_0} B^2 + \scro \left( \big( A^2 + |B| \big)^2 \right) 
  + \calO\left( \big( |c| + |\mu| \big) \big( |A| + |B| \big)^2 \right)
\]
as $(A,B,c,\mu) \to (0,0,0,0)$.

It remains to compute the expansion of the term
\[
- \rst{ \bigg( E_\mu(u) + \frac 1 2 s_\mu(u) K* s_\mu(u) \bigg) }{x=0}
\]
in \eqref{eq:NFE_Lyapunov}.
Since $\psi_{ijkl} \in \ker \calQ$, one has $\psi_{ijkl}(0) = \psi_{ijkl}'(0) = 0$, hence
\[
- \rst{ \bigg( E_\mu(u) + \frac 1 2 s_\mu(u) K* s_\mu(u) \bigg) }{x=0} = - E_\mu(A) - \frac 1 2 s_\mu(A) \rst{ K * s_\mu(u) }{x=0}.
\]
Since $u$ solves \eqref{eq:NFE_TW_modified}, one has
\begin{equation}
  \label{eq:computed_with_solution}
  \rst{ K*s_\mu(u) }{x=0} = \rst{ \big( u - c u' \big) }{x=0} = A - c B.
\end{equation}
With this one readily computes that
\[
  - \rst{ \bigg( E_\mu(u) + \frac 1 2 s_\mu(u) K* s_\mu(u) \bigg) }{x=0} = - \frac 1 4 \alpha A^4 \\
  + \calO\left(|A|^3 \big( |\mu| + |A| \big)^2 + |c||A||B| \right),
\]
as $(A,B,c,\mu) \to (0,0,0,0)$.

Combining these calculations, we obtain
\[
    \Lambda_{c,\mu}(A,B) = - \frac 1 4 \alpha A^4 - \frac 1 4 \frac{\kappa_2}{\kappa_0^2} B^2  + \scro \left( \big( A^2 + |B| \big)^2 \right) 
  + \calO\left( \big( |c| + |\mu| \big) \big( |A| + |B| \big)^2 \right),
\]
as $(A,B,c,\mu) \to (0,0,0,0)$.
Interestingly, \eqref{eq:computed_with_solution} is the only place in the computation of this expansion where we used that $\calM_0$ consists of solutions of \eqref{eq:NFE_TW_modified}.
Furthermore, we did not need to compute the Taylor expansion of $\Psi$ itself.

We now require that $(0,0)$ is a local maximum of $\Lambda_{c,\mu}$ by imposing the following restriction upon $K$.
\begin{myhyp}[Positive second moment]
  We demand that $\kappa_2 > 0$.
\end{myhyp}

Let $\epsilon , \epsilon' > 0$ be as dictated by the parameter dependent center manifold theorem, and fix $c \in (-\epsilon',\epsilon')\setmin\{0\}$.
Given $\delta' > 0$, define
\[
V := \set{ (A,B) \in \R^2 }{ \left| \frac 1 4 \alpha A^4 + \frac 1 4 \frac{\kappa_2}{\kappa_0^2} B^2 \right| < \delta' }.
\]
We choose $\delta'$ sufficiently small, so that
\[
\big\| u_0(x) + \Psi(u_0)(x) \big\|_{W^{1,2}_{-\delta}(\R,\R^d)} < \frac{\epsilon}{C_{0,2,-\eta}}, \qquad \text{for} \quad (A,B) \in V, \quad u_0 = A \e_0 + B \e_1,
\]
and the gradient of $\Lambda_{c,0}$ is inward pointing along $\bdy V$.
Here $C_{0,2,-\eta}$ is the constant from inequality \eqref{eq:Morrey_exponential}.
Since $\Lambda_{c,\mu}$ is monotone along the flow for $c \neq 0$,
we find that after decreasing $\delta'$ and $\epsilon'$ if needed, for any $\mu \in (-\epsilon',\epsilon')$ the flow $\chi_{c,\mu}^\tau$ is transverse to $\bdy V$.

\subsection{Whitham type equations}
\label{appendix:Whitham}

\subsubsection{Expanding the Hamiltonian}

\paragraph{Asymptotic expansion of $\Psi$.}
We formally expand $\Psi$ as
\[
\Psi(A,B,c-c_*) = \sum_{i,j,k \geq 1} A^i B^j (c-c_*)^k\psi_{i,j,k}, \qquad \text{where} \quad \psi_{i,j,k} \in \ker \calQ.
\]
The tangency $D_u \Psi(0,0,0) = 0$ implies $\psi_{1,0,0} = \psi_{0,1,0} = 0$.
Here we compute the terms of order $\calO\left( \big(|c-c_*||A| + A^2 \right)$ in the expansion of $\Psi$.
\begin{description}
\item[Terms $\calO(|c-c_*| |A|)$.] 
Collecting terms of order $\calO(|c-c_*| |A|)$, we find that
\[
c_* \psi_{1,0,1} + K* \psi_{1,0,1} = \e_0, \qquad \psi_{1,0,1} \in \ker \calQ.
\]
Substituting the ansatz
$\psi_{1,0,1}(x) = \beta_{1,0,1} x^2 + \varphi_{1,0,1}$ with $\varphi_{1,0,1} \in \calE_0$
and evaluating at $x=0$ leads to the compatibility condition $\beta_{1,0,1} \kappa_2 = 1$.
Hence, $\psi_{1,0,1}(x) = \kappa_2^{-1} x^2 + \varphi_{1,0,1}$.
Then $\calQ(\psi_{1,0,1}) = \varphi_{1,0,1}$, and ${\varphi_{1,0,1} = 0}$.

\item[Terms of order $\calO(A^2)$.]
At quadratic order in $A$, we find
\[
c_* \psi_{2,0,0} + K*\psi_{2,0,0} = - \alpha \e_0, \qquad \psi_{2,0,0} \in \ker \calQ.
\]
Substituting the ansatz $\psi_{2,0,0}(x) = \beta_{2,0,0} x^2 + \varphi_{2,0,0}$ with $\varphi_{2,0,0} \in \calE_0$
 and evaluating at $x=0$ leads to the compatibility condition
$\beta_{2,0,0} \kappa_2 = -\alpha$.
Hence $\psi_{2,0,0}(x) = -\alpha \kappa_2^{-1} x^2 + \varphi_{2,0,0}$.
Then $\calQ(\psi_{2,0,0}) = \varphi_{2,0,0}$, from which we find ${\varphi_{2,0,0} = 0}$.
\end{description}

\paragraph{Asymptotic expansion of $H_c$.}
We now compute the Taylor expansion of $H_{c}(A,B)$.

First, we consider the boundary term in \eqref{eq:Whitham_Hamiltonian}, that is,
\begin{equation}
  \label{eq:Whitham_boundary}
  \calB_c(u) = \frac 1 2 \iint_Q K(x-y)  \bigg( u(x) u'(y) - u'(x) u(y) \bigg) \d x \d y.
\end{equation}
Formally, we set
\[
\calB_c(u) = \sum_{i,j,k \geq 1} A^i B^j (c-c_*)^k b_{i,j,k}.
\]
In the following, we let $P$ and $\kappa_m$ be as defined in \eqref{eq:skew_pairing} and \eqref{eq:skew_pairing}, respectively.
We now compute the expansion of $\calB_\mu$.
\begin{description}
\item[Terms of order $\calO(A^2)$.]
At order $\calO(A^2)$ we obtain
\[
b_{2,0,0} = P(\e_0,\e_0') = 0.
\]
\item[Terms of order $\calO(B^2)$.]
At order $\calO(A^2)$ we find, using Lemma \ref{lemma:omega_covariance_computation},
\[
b_{0,2,0} = P(\e_1,\e_1') = \frac 1 2 C_{1,0} \kappa_2 = - \frac 1 4 \kappa_2.
\]
\item[Terms of order $\calO(|AB|)$.]
At order $\calO(|AB|)$ we find, using Lemma \ref{lemma:omega_covariance_computation} and $\kappa_1=0$ due to $K$ being even,
\[
b_{1,1,0} = \frac 1 2 P(\e_0,\e_1') = \frac 1 2C_{0,0} \kappa_1 = 0.
\]
\item[Terms of order $\calO(|A|^3)$.]
At order $\calO(|A|^3)$, we have
\[
b_{3,0,0} = P(\e_0,\psi_{2,0,0}') = - \frac{\alpha}{\kappa_2} C_{0,1} \kappa_2 = - \frac{\alpha}{2}.
\]
\item[Terms of order $\calO(A^2|B|)$.]
At order $\calO(A^2|B|)$, we find, using $\kappa_3=0$ due to $K$ being even,
\[
b_{2,1,0} = P(\psi_{2,0,0},\e_1') + P(\e_1,\psi_{2,0,0}') = \frac 1 2 C_{2,0} \kappa_3 - \frac{\alpha}{\kappa_2} C_{1,1} \kappa_3 = 0.
\]
\item[Terms of order $\calO(|c-c_*|A^2)$.]
Finally, the first nonzero parameter-dependent term is
\[
b_{1,0,1} = P(\e_0,\psi_{1,0,1}') = \frac{1}{\kappa_2} C_{0,1} \kappa_2 = \frac 1 2.
\]
\end{description}
We conclude that
\begin{multline*}
    \calB_c(u) = - \frac{\alpha}{2} A^3 + \frac 1 2 (c-c_*) A^2 - \frac 1 4 \kappa_2 B^2 \\
+\calO\left( A^4 + |B|^3 + |A| B^2+ |A|^3 |B|+|c-c_*|( |A|^3 + |AB| + B^2 )+ |c-c_*|^2 A^2  \right)
\end{multline*}
as $(A,B,c) \to (0,0,c_*)$.
It remains to compute the expansion of the term
\[
- \rst{ \bigg( F_c(u) + \frac 1 2 u K * u \bigg) }{ x=0 } 
\]
in \eqref{eq:Whitham_Hamiltonian}.
Since $\psi_{ijk} \in \ker \calQ$, one has $\psi_{ijk}(0) = \psi_{ijk}'(0) = 0$, hence
\[
- \rst{ \bigg( F_c(u) + \frac 1 2 u K * u \bigg) }{ x=0 } = - F_c(A) - \frac 1 2 A \rst{ K*u }{x=0}.
\]
Using the fact that $u$ solves \eqref{eq:Whitham_TW}, we obtain
\[
\rst{ K*u }{x=0} = - F'_{c}(A)  = c A - \alpha A^2 + \calO\left( |A| \big( |A| + |c-c_*| \big)^2 \right), \qquad \text{as} \quad (A,c) \to (0,c_*),
\]
hence
\[
- \rst{ \bigg(F_{c}(u) + \frac 1 2 u K * u \bigg) }{ x=0 }  = - F_{c}(A) + \frac 1 2 A F'_{c}(A) = \frac1 6 \alpha A^3 + \calO\left( A^2 \big( |A| + |c-c_*| \big)^2 \right)
\]
as $(A,c) \to (0,c_*)$.

Combining these calculations, we have obtained the expansion
\begin{multline*}
  H_c(A,B) = \frac{1-3\alpha}{6} A^3 + \frac 1 2 (c-c_*) A^2 - \frac 1 4 \kappa_2 B^2  \\
+\calO\left( A^4 + |B|^3 + |A| B^2+ |A|^3 |B|+|c-c_*|( |A|^3 + |AB| + B^2 )+ |c-c_*|^2 A^2  \right)
\end{multline*}
as $(A,B,c) \to (0,0,c_*)$.





\bibliography{biblio}
\bibliographystyle{abbrv}

\end{document}